\documentclass[11pt]{amsart}
 \usepackage{amsopn}
 \usepackage{amsmath,amsthm,amssymb}

\begin{document}

\large \pagestyle{plain}

\vspace{45mm}

\begin{center}

{\bf PROJECTIVE STRUCTURES AND EXACT VARIATIONAL FORMULA OF
MONODROMY GROUP OF THE LINEAR DIFFERENTIAL EQUATIONS ON COMPACT
RIEMANN SURFACE.}

\vspace{5mm}

{\bf V.V. Chueshev}

RUSSIA, 650043, Kemerovo, Kemerovo State University, Krasnaia, 6

e-mail: chueshev@lanserv1.kemsu.ru

\end{center}

\vspace{5mm}

\begin{abstract}

{\bf Abstract:} In this paper we are investigated the monodromy
group for linearly polymorphic functions on compact Riemann
surface of genus  $g \geq 2,$ in connection with standard
uniformization of these surfaces by Kleinian groups, and are found
a neccessary and sufficients conditions, that a linearly
polymorphic function on compact Riemann surface gave a standard
uniformization of this surface. We are investigated the monodromy
mapping $p : {\bf T}_{g}Q \rightarrow \mathcal{M},$ where ${\bf
T}_{g}Q$ is a vector bundle of holomorphic quadratic abelian
differentials over the Teichmueller space of compact Riemann
surfaces of genus $g,$ $\mathcal{M}$ is a space of monodromy
groups for of genus $g.$ Here is proved that over any space, which
consist from quasiconformal deformations by Koebe group of
signature $\sigma = (h, s; i_{1},..., i_{m}),$ connected with
standard uniformization compact Riemann surface of genus $g =
|\sigma|,$ this mapping $p$ has the lifting of path property .
Also we are received exact variational formula for monodromy group
of the linear differential equation of the second order and the
first variation for solution of the Schwartz equation on compact
Riemann surface.

\end{abstract}

{\bf AMS Subj.Classification:} 30F10, 30F30, 30F35, 30F40, 32G15


\vspace{3mm}

The monodromy groups of linearly polymorphic functions on compact
Riemann surface have appeared of the 19-th century in works by F.
Klein, H. Poincare, E. Picard [see 1-3], in connection with a
uniformization problem compact Riemann surface. In 70 years of
20-th century monodromy groups have appeared again in works by C.
Earle [4], I. Кra [5; 6], B. Maskit [7], D.A. Hejhal [1-3] and R.
Gunning [8], in connection with general uniformization problem and
with the theory Teichmueller spaces. In 1985 year P.G. Zograf,
L.А. Таchtajain [9] solved the problem of the accessory parameters
of the linear differential equation of the second order (so-called
Fuchsian type) on compact Riemann surface. А.B. Venkov [10] has
found the explicit formulas for these parameters in terms
monodromy groups, which are a Fuchsian group.

\vspace{3mm}

\centerline {{\bf 1. Uniformization and linearly polymorphic
functions.}}

\vspace{3mm}

In this section 1 we are investigated the monodromy group for
linearly polymorphic functions on compact Riemann surface of genus
$g \geq 2,$ in connection with standard uniformization of these
surfaces by Kleinian groups. Uniformization $(\Delta, G)$ for
compact Riemann surface $F$ is called to be  standard, if group
$G$ have not an elliptic elements and have not accidental
parabolic elements, or natural projection $\pi : \Delta
\rightarrow F$ is planar regular not ramified covering [7]. Here
we are shown, how the algebraic description of monodromy group is
connected with standard uniformization of compact Riemann surface.
Also we find a neccessary and sufficients conditions, that a
linearly polymorphic function on compact Riemann surface gave a
standard uniformization of this surface. These conditions have
simple topological sense.

Let $F$ be compact Riemann surface of genus $g \geq 2,$
$\pi_{1}(F, O)$ - fundamental group surface $F$ with basic point
$O.$ A complex projective structure on $F$ is a maximal atlas from
cards on $F$ such that all transition maps belong to the group
$PSL(2, {\bf C})$ [8].

{\bf Definition 1.1.} A locally meromorphic (multivalued) function
$z$ on $F,$ which transforms linear fractionally under action
group $\pi_{1}(F, O),$ is called linearly polymorphic functions on
$F.$

Further we are considered only locally schlicht linearly
polymorphic function on $F.$ They play a very important role in
modern uniformization theory and in theory of the Teichmueller
spaces.

Let $(U, \pi)$ be the universal covering for $F,$ where $\pi : U
\rightarrow F$ - natural projection, $U$ - unit disc with the
center in the origin of coordinates on the extended complex plane
$\overline{\bf C},$ and $\Gamma$ is a Fuchsian group such that $F
= U/\Gamma .$

The linearly polymorphic function $z$ on $F$ can be lifted to a
meromorphic (single-valued) function $z = z(t)$ on $U$ in such a
way  that
$$
z(A t) = \widetilde{A} z(t), t \in U, A \in \Gamma,
$$
where $\rho(A) = \widetilde{A} \in PSL(2, {\bf C})$ is group of
linear fractionally transformations. Therefore a monodromy
homomorphism $\rho : \Gamma \rightarrow PSL(2, {\bf C})$ is
determined. The image $\rho(\Gamma)$ for $\Gamma$ on $\rho$ is
called  monodromy group $\mathcal{M}[z]$ for $z$ on $F.$ I. Kra
[6] has named this function to be $(\Gamma, \rho)-$deformation of
Fuchsian group $\Gamma$ in $U,$ аnd R. Gunning [8](see also М.E.
Каpovich [11]) is called it a developing map of this unramified
complex projective structure on $F.$

For the theory of functions more suitable is the term linearly
polymorphic function, accepted in work D.А. Hejhal. This is
according to historical tradition, coming from H. Poincare, P.
Аppell and E. Coursat [12].

Let us give two classical examples of linearly polymorphic
function, which are connected with uniformization compact Riemann
surfaces.

Let $\pi : U \rightarrow F$ be universal covering mapping, for
which $\Gamma$ is group of covering transformations. Тhen
multivalued inverse function  $w = {\pi}^{-1}(\xi)$ will be
locally schlicht linearly polymorphic function on $F$ with
monodromy group $\Gamma.$

Let $D \subset \overline{\bf C}$ and $ \pi_{S} : D \rightarrow F$
is a Schottky covering, where the group of covering
transformations of $\pi_{S}$ is the Schottky group $\Gamma_{S}.$
Тhen function $w = \pi_{S}^{-1}(\xi)$ will be locally
single-valued linearly polymorphic function on $F$ with the
monodromy group $\Gamma_{S}.$

Let us note, that there are  monodromy groups, which algebraically
are arranged as a Fuchsian group $\Gamma$ or as a Schottky group
$\Gamma_{S},$ but not connected with uniformization compact
Riemann surface [1].

The function
$$
q(t) = \{z, t\} \equiv (\frac{z''(t)}{z'(t)})' -
\frac{1}{2}(\frac{z''(t)}{z'(t)})^{2}
$$
is a  Schwartz derivative for function $z = z(t)$ on $U.$ This
function is  holomorphic on $U$ and satisfies the condition
$$
q(A t) (d A t)^{2} = q(t) dt^{2}, t \in U, A \in \Gamma.
$$
Hence, $q = q(t) dt^{2}$ is the holomorphic (abelian) quadratic
differential on $U/\Gamma = F.$

Such linearly polymorphic functions on $F$ can be received, as the
quotient of two the linearly independent solutions of the linear
differential equation of the second order Fuchsian type on $F$
[1].

Let  $B_{2}(U, \Gamma)$ be a complex vector space of holomorphic
quadratic differentials  $q = q(t) dt^{2}$ with respect to
$\Gamma$ on $U$ with norm
$$
\| q \| = sup _{t \in U} \{(1 - |t^{2}|)^{2} | q(t) | \} < \infty.
$$
Classical method to construct of linearly polymorphic function on
quadratic differential consist in the following: let $q \in
B_{2}(U, \Gamma),$ then there exists the unique locally schlicht
meromorphic function $z = z(t)$ on $U$ such that $\{z, t\} = q(t)$
and $z(t) = \frac{1}{t} + O(|t|)$ in neighborhood point $t = 0.$
Hence for any $A \in \Gamma$ there is unique element $\rho_{z}(A)
\in PSL(2, {\bf C})$ such that $z(A t) = \rho_{z}(A) z(t), t \in
U$ [1].

Let $\tau_0 = [F, \{a_{k}, b_{k}\}^{g}_{k=1}]$ be a  marked
compact Riemann surface of genus $g \geq 2.$ We will choose and
fixed a point $t_{0} \in U$ and set $O = \pi(t_{0}).$ Then there
exists a natural isomorphism between $\pi_{1}(F, O)$ and marked
Fuchsian group
$$
\Gamma = \{A_{1},..., A_{g}, B_{1},..., B_{g} : [A_{1},
B_{1}]...[A_{g}, B_{g}] = 1 \}
$$
of the first kind on $U,$ defined by $a_{j}\rightarrow A_{j},
b_{j}\rightarrow B_{j}, j = 1,..., g.$

From the theory of Teichmueller spaces ${\bf T}_{g}$ it is known
that there is a homeomorphism translating $\tau = [F_{\tau},
\{a_{k}(\tau), b_{k}(\tau)\}^{g}_{k=1}] \in {\bf T}_{g}$ in group
$$
\Gamma_{\tau} = \{A_{1}(\tau),..., B_{g}(\tau) :
 \prod_{j=1}^{g}[A_{j}(\tau), B_{j}(\tau)] = 1 \}
$$
from space of normalized marked Fuchsian groups on $U$ [13]. Let
$z$ be a linearly polymorphic function on $F_{\tau} =
U/\Gamma_{\tau},$ then for meromorphic function $z = z(t)$ on $U$
we have $z(A t) = \widetilde{A} z(t), A \in \Gamma_{\tau},$
$\widetilde{A} \in PSL(2, {\bf C}).$ The mapping $A \rightarrow
\widetilde{A}$  is called the monodromy homomorphism. It
determines the marked monodromy group
$$ \mathcal{M}[z] =
\{\widetilde{A}_{1}(\tau),..., \widetilde{B}_{g}(\tau) :
[\widetilde{A}_{1}(\tau),
\widetilde{B}_{1}(\tau)]...[\widetilde{A}_{g}(\tau),
\widetilde{B}_{g}(\tau)] = 1\},
$$
i.e. $\mathcal{M}[z]$ is a point in $[PSL(2, {\bf C})]^{2g}.$

Let $\mathcal{F}_{\tau}$ be a fundamental polygon for group
$\Gamma_{\tau}$ in $U,$ whose border
$$
\partial \mathcal{F}_{\tau} = a_{1}^{+}(\tau) b_{1}^{+}(\tau)
a_{1}^{-}(\tau) b_{1}^{-}(\tau)... a_{g}^{+}(\tau) b_{g}^{+}(\tau)
a_{g}^{-}(\tau) b_{g}^{-}(\tau)
$$
is lifting of commutator path
$$
[a_{1}(\tau), b_{1}(\tau)]...[a_{g}(\tau), b_{g}(\tau)]
$$
from $t_{0} \in U;$ the sides of $\mathcal{F}_{\tau}$ are paired
identified by transformations $A_{k}(\tau) : a_{k}^{-}(\tau)
\rightarrow a_{k}^{+}(\tau), B_{k}(\tau) : b_{k}^{-}(\tau)
\rightarrow b_{k}^{+}(\tau), k = 1,..., g.$

According to F.Klein, we define the fundamental membrane $R_{z}$
for function $z = z(t)$ as the Riemannian (multivalued) image
$z(\mathcal{F}_{\tau}).$ Its is simply connected and unramified,
but it is possible for mapping $z : \mathcal{F}_{\tau} \rightarrow
z(\mathcal{F}_{\tau})$ to be $n-$valued, $n \geq 2.$ The sides of
membrane are paired identified by transformations
$\widetilde{A}_{1}(\tau),..., \widetilde{B}_{g}(\tau).$ There is
an equivalence, connecting linearly polymorphic function $z$ and
its fundamental membrane. So, given a simply connected unramified
domain $R$ with indicated properties, there exists a locally
schlicht linearly polymorphic function $z$ on some marked compact
Riemann surface $F_{\tau}$ such that $R_{z} = R$ [1].

Explicit construction of such a function $z$ on $R$ is indicated
in [1], and there is its generalizations, where instead of $(U,
F_{\tau})$ it one take the standard uniformization $(\Delta, G)$
for compact Riemann surface of genus $g.$ Indeed, let $z = z(t)$
be a locally schlicht linearly polymorphic function on $\tau \in
{\bf T}_{g}$ such that
$$
\mathcal{M}[z] = \{ T_{1},..., T_{h}, 1,..., 1, U_{1},...,
U_{g-h},
 V_{1},..., V_{g-h} :  [U_{1}, V_{1}] = ... = [U_{s}, V_{s}] =
$$
$$
 \prod_{j=1}^{i_{1}}[U_{s+j}, V_{s+j}] = ...
 = \prod_{j=1}^{i_{m}}[U_{g-i_{m}+j}, V_{g-i_{m}+j}] = 1 \}.
\eqno(*)
$$
Suppose that $w : U \rightarrow D_{\tau}$ is a  Koebe
uniformization of signature $\sigma = (h, s; i_{1},..., i_{m})
\neq (0, 2; 0, ..., 0),$ $|\sigma| = h + s + i_{1} + ... + i_{m} =
g,$ $i_{j} \neq 1, j = 1,..., m,$ for marked compact Riemann
surface $\tau$ and
$$
G_{\tau} = \{T_{1}',..., T_{h}', U_{1}',..., U_{g-h}',
 V_{1}',..., V_{g-h}'
 : [U_{1}', V_{1}'] = ... = [U_{s}', V_{s}'] = ...
$$
$$
 = \prod_{j=1}^{i_{1}}[U_{s+j}', V_{s+j}'] = ...
 = \prod_{j=1}^{i_{m}}[U_{g-i_{m}+j}', V_{g-i_{m}+j}'] = 1\}
$$
is the corresponding marked Koebe group of signature $\sigma.$
Then the function $Z = z w^{-1}$ is locally schlicht linearly
polymorphic function on $D_{\tau}$ (invariant component of group
$G_{\tau}).$ It satisfies the following relations:
$$
Z T_{k}'(\widetilde{t}) = T_{k} Z(\widetilde{t}), k = 1,..., h, Z
U_{j}'(\widetilde{t}) = U_{j} Z(\widetilde{t}),
$$
$$
Z V_{j}'(\widetilde{t}) = V_{j} Z(\widetilde{t}), j = 1,..., g -
h.
$$
Let $K_{\tau}$ be a standard fundamental $(2h + s + m)-$connected
domain for group $G_{\tau}$ which sides are identified (in pairs)
by standard generators of group $G_{\tau},$ changing orientation
of determining curves.

We list the basic properties of a fundamental membrane
$Z(K_{\tau}):$

1) $Z(K_{\tau})$ is planar (i.е. there is conformal mapping its on
planar domain) and $(2h + s + m)-$connected,

2) $Z(K_{\tau})$ is unramified,

3) the sides of $Z(K_{\tau})$ are identified, by linear
fractionally transformation
$$
T_{1},..., T_{h}, 1,..., 1, U_{1},..., U_{g-h},
 V_{1},..., V_{g-h},
$$
with relations (*). Repeating topological and analytical
construction by D.А. Hejhal [1, c. 28-29], we have receive that
for any membrane $R$ with properties 1) - 3) exists linearly
polymorphic functions $z = z(t)$ on $U$ and $Z(\widetilde{t})$ on
$D_{\tau}$ for some $\tau \in {\bf T}_{g},$ which satisfy to the
equation
$$
z = Z w,
$$
where $z(U) = Z(D_{\tau}).$ Hence, we have an equivalence
$$
(z(t), Z(\widetilde{t})) \longleftrightarrow R,
$$
where $R$ has properties 1) - 3).

{\bf Theorem} (Kra - Gunning) [8; 5; 6].{ \sl{ Let $z$ be a
locally schlicht linearly polymorphic function on compact Riemann
surface $U/\Gamma = F$ of genus $g \geq 2.$ Then the following
conditions are equivalent:

1) $\mathcal{M}[z]$ acts on $z(U)$  discontinuously,

2) $z : U \rightarrow z(U)$ is a (topological) covering mapping,

3) $z(U) \neq \overline{\bf C}.$}}

{\bf Lemma 1.1.} {\sl{ Let $w = w(t)$ be a  locally schlicht
linearly polymorphic function on compact Riemann surface $F$ of
genus $g \geq 2$ and $w(U) \neq \overline{\bf C}.$ Then
$\mathcal{M} [w]$ is nonelementary, finite generated Kleinian
group with invariant component $w(U).$}}

{\sl Proof.} Nonelementarity follows from theorem 7 [8] since
$w(U)$ has a hyperbolic covering. Discontinuity of
$\mathcal{M}[z]$ or that $\mathcal{M}[z]$ is Kleinian group,
follows from the theorem Kra-Gunning. As a consequence from this
theorem we have that $\mathcal{M}[z]$ can not act discontinuously
in greater domain than $w(U),$ i.e. the domain $w(U)$ is invariant
component for group $\mathcal{M}[z].$ Lemma is proved.

The main result in this section 1 is the following

{\bf Theorem 1.2.} {\sl{Let $w = w(t)$ be a locally schlicht
linearly polymorphic function on  compact Riemann surface $F$ of
genus $g \geq 2.$ Then $w = w(t)$ is an uniformization of $F$ if
and only if the following conditions are carried out:

1) $w(U) \neq \overline {\bf C},$

2) $w(U)/ \mathcal{M}[w]$ is a compact surface of genus $g.$}}

{\sl Proof.} If $w = w(t)$ be an uniformization of $F,$ then, by
definition, pair $(w(U), \mathcal{M}[w])$ is those that $w(U)/
\mathcal{M}[w]$ is conformal equivalent to $F,$ and $w(U)$ is the
covering surface for $F.$ Hence, universal covering surface for
$w(U)$ will be a disc and so $w(U) \neq \overline {\bf C}.$

Conversely, from the condition $w(U) \neq \overline {\bf C}$ it
follows that $w : U \rightarrow w(U)$ is a  topological covering.
Consider the commutative diagram of mapping given below. Here
$\pi$ и $\widetilde{\pi}$ are natural projections and $\widetilde
w$ is proper holomorphic covering mapping from $F$ on compact
Riemann surface $F_{1}$ [16]. We notice that $\widetilde{\pi}$ is
ramified if and only if  $\widetilde w$ is ramified.

\centerline {\bf }

\begin{center}

\begin{picture}(140,140)
\put(85,32){\makebox(0,0)[tr]{$\widetilde{w}$}}
\put(92,130){\makebox(0,0)[tr]{$w$}}
\put(-10,95){\makebox(0,0)[tr]{$\pi$}}
\put(215,95){\makebox(0,0)[tr]{$\widetilde{\pi}$}}

\put(20,45){\makebox(0,0)[tr]{$U/\Gamma = F$}}

\put(250,45){\makebox(0,0)[tr]{$F_{1} = w(U)/\mathcal{M}[w]$}}

\put(0,140){\makebox(0,0)[tr]{$U$}}

\put(220,140){\makebox(0,0)[tr]{$ w(U)$}}

\put(200,120){\vector(0,-1){60}}

\put(-5,120){\vector(0,-1){60}}

\put(30,135){\vector(1,0){130}}

\put(40,40){\vector(1,0){90}}

\end{picture}
\end{center}

From the second condition of the theorem and Riemann-Hurwitz
formula [17; 16] we are received that $\widetilde w$ be one-to-one
and $\widetilde w$ is a conformal homeomorphism $F$ on $F_{1}.$
Theorem is proved.

Remark 1.1. It is possible to prove that from  condition 2) of the
theorem 1.2 implies the condition 1). Inverse is not correct, see
for example [6].

From the theorem 1.2 and the theorem 3 [22] we obtain following

{\bf Corollary 1.3.} {\sl{ Let $w = w(t)$ is a locally schlicht
linearly polymorphic function on marked compact Riemann surface
$F$ of genus $g \geq 2$ such that $\mathcal{M} [w]$ is a marked
Koebe group of signature $\sigma = (h, s; i_{1}, ..., i_{m}) \neq
(0, 2; 0, ..., 0),$ $|\sigma| = g,$ $i_{j} \neq 1, j = 1,..., m.$
If $w(U) \neq \overline{\bf C},$ then $w = w(t)$ is a
uniformization of $F$ by the Koebe group of signature $\sigma.$}}

Remark 1.2. If $w = w(t)$ is such as in the corollary 1.3, but
$w(U) = \overline{\bf C},$ then statement of the corollary 1.3 is
not truly. For the proof it is enough to construct a fundamental
membrane $R_w$ for $w$ such that $\mathcal{M}[w]$ is the same as
in  corollary 1.3, but the mapping $w: \mathcal F \rightarrow R_w
= w(\mathcal F)$ is $n$-valued, $n \geq 2.$ Here $\mathcal F$ is
an fundamental polygon of Fuchsian group $\Gamma_{w},$ which
uniformize $F$ in disc $U.$ For signature $\sigma = (h, s;
i_{1},...,i_{m}) \neq (0, 2; 0,..., 0),$ $|\sigma| = g \geq 2,$ we
will consider the following cases:

1) if $h \geq 1,$ then $R_w$ can be construct as  in the theorem 6
[1] in the form of two-sheeted covering of a ring on ambient
surface $\{w^2 = z \},$ and from the top sheet it is necessary to
remove domains which are bounded by: $2(h - 1)$ by the closed
curves, $s$ curvilinear quadrangles and $m$ curvilinear polygons
with $4i_1,..., 4i_m$ sides;

2) if $h = 0, s \geq 2,$ then $R_w$ can be construct as
two-sheeted ring domains, which are bounded two curvilinear
quadrangles, and further - as in 1);

3) case $h = 0, s = 1, i_1 = ... = i_m = 0$ does not meet in to
kind that $g \geq 2;$

4) if $h = 0, s = 0, 1,$ then exists $k$ such that $i_k \geq 2,$
and $R_w$ can be construct as in the theorem 4 [1], as two-sheeted
covering on ambient surface $\{w^2 = (z - c_1) (z - c_2), c_1, c_2
\in {\bf C}, c_1 \neq c_2\},$ and further - as in 1).

Remark 1.3. The corollary 1.3 and remark 1.2 show, that, as well
as in classical problem of a choise of the accessory parameters
for uniformization by Fuchsian groups, problem of a choise of
accessory parameters [1] for any standard uniformization compact
Riemann surface by Koebe groups has unique solution, up to linear
fractionally transformation, if  linearly polymorphic function $w$
has the limited image of a disc, i.e. $w(U) \neq \overline {\bf
C}.$

Remark 1.4. The corollary 1.3 include, in particular, the theorem
3 and 5 by D.A. Hejhal [1] for the Fuchsian group and the Schottky
group respectively. Under the remark 1.2 we receive examples of
monodromy groups, which are algebraically arranged as the marked
Koebe group, but they are act non discontinuously on $w(U)$,
though $\mathcal{M}[w]$ is group of conformal homeomorphisms
many-sheeted Riemannian domain $w(U)$ on itself.

\vspace{3mm}

\centerline{{\bf 2. The monodromy mapping}}

\vspace{3mm}

In this section we are investigated in the monodromy mapping $p :
{\bf T}_{g}Q \rightarrow \mathcal{M},$ where ${\bf T}_{g}Q$ is a
vector bundle of holomorphic quadratic abelian differentials over
the Teichmueller space of compact Riemann surfaces of genus $g,$
and $\mathcal{M}$ is a space of monodromy groups for of genus $g.$
D.A. Hejhal in [1] has shown that mapping $p$ is local
homeomorphism, but has not by the lifting of path property over
$\mathcal{M}.$ But over a part $\mathcal{M}_{q},$ consisting of
quasifuchsian uniformizations, it already has this property.
Naturally it is interesting to find parts of space $\mathcal{M},$
admitting this property. In this section is proved that over any
space, which consist of quasiconformal deformations by Koebe group
of signature $\sigma = (h, s; i_{1},..., i_{m}),$ connected with
standard uniformization compact Riemann surface of genus $g =
|\sigma|,$ this mapping $p$ has the lifting of path property.

Two ordered collections  $X_{1},..., X_{g}, Y_{1},..., Y_{g}$ and
$X'_{1},..., X'_{g}, Y'_{1},..., Y'_{g}$ are $PSL(2, {\bf
C})-$equivalent, if exists $B \in PSL(2, {\bf C})$ such that
$X'_{k} = B X_{k} B^{-1}, Y'_{k} = B Y_{k} B^{-1}$ for all $k =
1,..., g.$ Let $\widetilde{\mathcal{M}}$ be the set of marked
monodromy groups $\mathcal{M}[z]$ for all $\tau \in {\bf T}_{g}.$
Set $\mathcal{M} = \widetilde{\mathcal{M}} mod PSL(2, {\bf C}).$
This space we will named as space of the marked monodromy groups
the fixed of genus $g.$ In $\widetilde{\mathcal{M}}$ and
$\mathcal{M}$ can be defined the topology of convergence on
generators and the factor-topology respectively. R. Gunning [8]
and D.A. Hejhal [1] have proved that $\mathcal{M} \subset
\mathcal{N} \subset [PSL(2, {\bf C})]^{2g} mod PSL(2, {\bf C}),$
where $\mathcal{N}$ be complex-analytic manifold of complex
dimension $6g - 6$ and $\mathcal{N}$ be Hausdorf locally
metrizable space and $\mathcal{M}$ be subdomain in $\mathcal{N}.$

For $\tau \in {\bf T}_{g}$ we will denote through $Q(\tau)$ the
complex vector space of holomorphic quadratic differentials $q =
q(t,\tau) dt^{2}$ on $U/\Gamma_{\tau}.$ The solution $z = z(t)$ of
the Schwartz equation $\{z, t\} = q(t, \tau)$ is determined up to
transformation from $PSL(2, {\bf C}).$ Such a way, the monodromy
mapping is well defined by
$$
p(\tau, q(t,\tau) d t^{2}) = \mathcal{M}[z] mod PSL(2, {\bf C}).
$$

It is well known that: 1) $dim_{\bf C} Q(\tau) = 3g - 3,$ 2) it is
possible to enter basis $q_{k}(t,\tau) dt^{2}, k = 1,..., 3g - 3,$
in $Q(\tau),$ which is globally complex-analytic depending from
$\tau$ on ${\bf T}_{g}$[18]. Let $ {\bf T}_{g}Q$ be holomorphic
vector bundle of holomorphic quadratic differentials over
complex-analytic manifold ${\bf T}_{g}.$ Hence, the monodromy
mapping is determined
$$
p : {\bf T}_{g}Q \rightarrow \mathcal{M}.
$$

It is easy to see that $dim_{\bf C}{\bf T}_{g}Q = 6g - 6.$ Indeed
by the theorem of H. Grauert [19], by virtue of simply
connectivity of ${\bf T}_{g},$ the bundle is analytically
equivalent to the trivial vector bundle of rank $3g - 3$ over
${\bf T}_{g}.$

A Beltrami differential with respect to a Kleinian group $G$ is
the form $\mu = \mu(z) d\overline{z}/d z,$ where

1) $\mu(z) \in L_{\infty}(\bf C),$

2) $\mu(A z)\overline{A'(z)}/A'(z) = \mu(z), A \in G,$

3) $\mu|\Lambda(G) = 0,$ where $\Lambda(G)$ be limit set of group
$G.$ We will denote through $M(\Delta, G) = \{\mu : supp \mu
\subset \Delta \}$ a complex Banach space with the norm $\|\mu
\|_{\infty}$ and $M_{0}(\Delta, G) = \{ \mu \in M(\Delta, G) :
\|\mu \|_{\infty} < 1 \}.$ Set that $0, 1, \infty \in \Lambda(G).$
We will consider a quasiconformal automorphisms $w = f^{\mu}(z)$
on plane, where $\mu \in M_{0}(\Omega(G), G),$ i.e. quasiconformal
deformations of group $G.$ Every automorphisms $f^{\mu}$ generate
a Kleinian group $G_{\mu} = f^{\mu} G (f^{\mu})^{-1}$ and
isomorphism $\chi_{\mu} : G \rightarrow G_{\mu}.$ Such two
automorphisms $f^{\mu_{1}}$ and $f^{\mu_{2}},$ and also Beltrami
differentials $\mu_{1}$ and $\mu_{2}$ appropriate to them from
$M_{0}(\Omega(G), G),$ we will name (strongly) quasiconformal
equivalent, if

1) they are homotopic on each surface $F_{j} \subset \Omega(G)/G,$

2) $f^{\mu_{1}} = f^{\mu_{2}}$ on $\Lambda(G)$ [20].

It is clear that under these conditions $G_{\mu_{1}} =
G_{\mu_{2}}$ and $\chi_{\mu_{1}} = \chi_{\mu_{2}}.$ Set of classes
of quasiconformal equivalent Beltrami differentials $\mu \in
M_{0}(\Omega(G), G)$ is called the space $\widehat{\bf
T}(\Omega(G), G)$ of quasiconformal deformations of group $G.$ Let
$\Delta$ be an $G-$invariant union connected components from
$\Omega(G).$ We will name  $\mu_{1}$ and $\mu_{2}$ from
$M_{0}(\Delta, G)$ to be of weakly quasiconformal equivalent, if
$f^{\mu_{1}} = f^{\mu_{2}}$ on $\Lambda(G)$ (after suitable linear
fractionally normalization). Factor-space ${\bf T}(\Delta, G)$
space $M_{0}(\Delta, G)$ under this relation of equivalence we
will name the space of weak quasiconformal deformations of group
$G$ with supports in $\Delta.$ Note that if all components from
$\Delta$ will be simply connected, then ${\bf T}(\Delta, G) =
\widehat{{\bf T}}(\Delta, G)$ [20].

Denote through ${\bf T}(G_{\sigma}) \equiv {\bf
T}(\Omega(G_{\sigma}), G_{\sigma})$ the space of weak
quasiconformal deformations of marked Koebe group $G_{\sigma}$ of
signature $\sigma,$ which give standard uniformization in
invariant component $\Delta_{\sigma}$ of marked compact Riemann
surface $F_{\sigma}$ of genus $g \geq 2,$ where $|\sigma| = g.$

{\bf Definition 2.1.} Appropriate to pair $(\tau, q) \in {\bf
T}_{g}Q$ the solution $z : U \rightarrow \overline{\bf C}$ of the
Schwartz equation $\{z, t\} = q(t, \tau)$ is defined normalized
quasiconformal deformation $f,$ $[\mu_{f}] \in {\bf
T}(G_{\sigma}),$ if exists  $A \in PSL(2, {\bf C})$ such that:

1) $z(\mathcal{F}_{\tau}) = A f(\mathcal{F}_{\Delta_{\sigma}}),$
where $\mathcal{F}_{\tau}$ and $\mathcal{F}_{\Delta_{\sigma}}$ be
standard fundamental domains for uniformizations for $\tau$ and
$F_{\sigma}$ by Fuchsian group $\Gamma_{\tau}$ in $U$ and by
$G_{\sigma}$ in $\Delta_{\sigma}$ respectively,

2) $\mathcal{M}[z] = A f G_{\sigma} f^{-1} A^{-1}$ (equality the
marked groups).

This definition is correct as if exists $f_{1}, f_{2}$ appropriate
to pair $(\tau, q)$ and $z,$ then from 2) follows $A_{1} f_{1}
\equiv A_{2} f_{2}$ on $\Lambda(G_{\sigma})$ for some $A_{1},
A_{2} \in PSL(2, {\bf C}),$ i.e. $\mu_{f_{1}} \in [\mu_{f_{2}}]$
in ${\bf T}(G_{\sigma}).$

We will identify ${\bf T}(G_{\sigma})$ with his image in
$\mathcal{M}$ on mapping
$$
{\bf T}(G_{\sigma}) \ni [\mu] \rightarrow (\chi_{\mu}(A'_{1}),
\chi_{\mu}(B'_{1}),..., \chi_{\mu}(A'_{g}), \chi_{\mu}(B'_{g}))
\in \mathcal{M},
$$
where $G_{\sigma} = \{A'_{1}, B'_{1},..., A'_{g}, B'_{g}\}$ be
fixed marked Koebe group $G_{\sigma}$ of signature $\sigma = (h,
s; i_{1},..., i_{m}),$ $g = |\sigma|,$ and its has $h$ generators
equal 1, for example, $A'_{1} = ... = A'_{h} = 1.$ Without
restriction of a generality we assume that $0, 1, \infty \in
\Lambda(G_{\sigma}).$ Here $\chi_{\mu}(A) = w^{\mu} A
(w^{\mu})^{-1},$ $A \in G_{\sigma},$ $ \mu$ is the Beltrami
differential for $G_{\sigma}$ on $\Omega(G_{\sigma})$ and
$w^{\mu}$ is a quasiconformal automorphism on $\overline{\bf C},$
$w^{\mu}(0) = 0,$ $w^{\mu}(1) = 1,$ $w^{\mu}(\infty) = \infty,$
being the solution of the Schwartz equation with coefficient $\mu$
[20].

We note that for the $\sigma = (h, s; i_{1},..., i_{m})$ the space
${\bf T}(G_{\sigma})$ is complex-analytic submanifold of complex
dimension $3g - 3 + 3(i_{1} + ... + i_{m} - m)$ in
complex-analytic manifold $\mathcal{M}$ of complex dimension $6g -
6$ [20; 15].

Remind the following  classical results

{\bf Theorem} (H. Poincare [12]). {\sl{If $z = z(t)$ and $w =
w(t)$ are locally schlicht  linearly polymorphic function on the
same marked compact Riemann surface $\tau \in {\bf T}_{g},$ and
$\mathcal{M}[z] = \mathcal{M}[w],$ then $z(t) = w(t)$ on $U.$}}

The main results of this section is

 {\bf Theorem 2.1.} {\sl{The
monodromy mapping $p : {\bf T}_{g}Q \rightarrow \mathcal{M}$ has
the lifting of path property over any space ${\bf
T}(\Omega(G_{\sigma}), G_{\sigma}) \subset \mathcal{M},$ where the
marked Koebe group $G_{\sigma}$ standardly uniformize a marked
compact Riemann surface of genus $g \geq 2,$ in the invariant
component $\Delta_{\sigma},$ $\sigma = (h, s; i_{1}, ..., i_{m})
\neq (0, 2; 0, ..., 0),$ $i_{j} \neq 1, j = 1,..., m$ and
$|\sigma| = g.$}}

{\sl{Proof.}} Let a signature $\sigma \neq (0, 0; g),$ i.e.
$G_{\sigma}$ is not a Fuchsian group, since this case is
considered in [1]. Denote by ${\bf B}_{\sigma}$ the set of pairs
$(\tau, q(t,\tau)d t^{2}) \in {\bf T}_{g} Q$  such that any
solution of the Schwartz equation $\{z, t \} = q(t, \tau)$ is
defined by normalized quasiconformal deformation of group
$G_{\sigma}.$ By the theorem 1[1] the mapping $p : {\bf
B}_{\sigma} \rightarrow p({\bf B}_{\sigma})
 \equiv {\bf T}(G_{\sigma})$  is continuous and is a local homeomorphism.
But mapping $p : {\bf B}_{\sigma} \rightarrow {\bf T}(G_{\sigma})$
is not one-to-one for  $h \neq 0.$ Indeed, let  $(\tau_{1},
q_{1})$ and $(\tau_{2}, q_{2})$ be distinct elements from ${\bf
B}_{\sigma}$ such that $\mathcal{M}[z_{1}] = \mathcal{M}[z_{2}]$
(equality marked monodromy groups). Hence $z_{1}(U) \neq
\overline{\bf C}$ and $z_{2}(U) \neq \overline{\bf C}.$ By the
Maskit's theorem [21] there exists a conformal similarity of these
groups on the marked Koebe group $\widetilde{G}_{\sigma}$ same
signature $\sigma.$ Hence, under the theorem 3 [22] we have
$\tau_{1} = [F, \{a_{k}, b_{k}\}^{g}_{k=1}]$ and $\tau_{2} = [F,
\{\psi(a_{k}), \psi(b_{k}\}^{g}_{k=1}],$ where $F$ - compact
Riemann surface of genus $g,$ $\psi(O) = O$ and $\psi$ belongs to
not trivial group $\Theta_{\sigma}$ (see [22]) of homeomorphisms
surface $F$ on itself. If $\tau_{1} = \tau_{2}$ in ${\bf T}_{g},$
then under the theorem Poincare we will receive $q_{1} = q_{2}.$
From here $\tau_{1} \neq \tau_{2}$ in ${\bf T}_{g},$ i.e. $\psi$
is not homotopically identical mapping on $F$ and distinct points
$(\tau_{1}, q_{1}),$ $(\tau_{2}, q_{2})$ of ${\bf B}_{\sigma}$
under mapping  $p$ pass in the same point $\mathcal{M}[z_{1}] =
\mathcal{M}[z_{2}]$ in $\mathcal{M}.$

By the Poincare theorem the mapping $p$ on fiber ${\bf
B}_{\sigma}(\tau_{0})$ of set ${\bf B}_{\sigma}$ is homeomorphism
for any $\tau_{0} \in {\bf T}_{g}.$ It is  possible to identify
the fiber ${\bf B}_{\sigma}(\tau_{0})$ with space ${\bf
T}(\overline{\bf C} \backslash \Delta (G_{0,\sigma}),
G_{0,\sigma})$ by a weak quasiconformal deformations of marked
Koebe group $G_{0,\sigma},$ with supports on $\overline{\bf C}
\backslash \Delta(G_{0,\sigma}),$ where group $G_{0,\sigma}$ is of
signature $\sigma$ and standardly uniformize $\tau_{0}$ in
invariant component $\Delta(G_{0,\sigma}).$ Since, all non
invariant component of $G_{0,\sigma}$ are simply connected, we
have
$$
{\bf T}(\overline{\bf C} \backslash
\Delta(G_{0,\sigma}),G_{0,\sigma}) = \hat{{\bf T}}(\overline{\bf
C} \backslash \Delta(G_{0,\sigma}), G_{0,\sigma}),
$$
where $\hat{{\bf T}}(\overline{\bf C} \backslash
\Delta(G_{0,\sigma}), G_{0,\sigma})$ is the space of
quasiconformal deformations of group $G_{0,\sigma}$
 with supports in $\overline{\bf C} \backslash
\Delta(G_{0,\sigma})$ [20]. Therefore ${\bf T}(\overline{\bf C}
\backslash \Delta(G_{0,\sigma}), G_{0,\sigma}),$ as well ${\bf
B}_{\sigma}(\tau_{0}),$ is simply connected. Since fiber ${\bf
B}_{\sigma}(\tau_{0})$ are contractible, ${\bf B}_{\sigma}$ is
also contractible and simply connected[1].

Now we show that any continuous path $\gamma = \{\mathcal{M}(\xi)
: 0 \leq \xi \leq 1\}$ in  ${\bf T}(G_{\sigma})$ is  lifted to
continuous path $\widetilde{\gamma}$ in  ${\bf T}_{g}Q$ from any
point laying over $\mathcal{M}(0).$ Let $(\tau_{0}, q_{0}) \in
{\bf B}_{\sigma}$ lays over $\mathcal{M}(0) = \mathcal{M}[z_{0}],$
where $\{z_{0}, t\} = q_{0}(t, \tau_{0}).$ Denote by
$\Gamma_{\tau_{0}}$ a Fuchsian group, which uniformize $\tau_{0}$
with fundamental polygon $\mathcal{F}_{\tau_{0}}$ in disc $U.$ We
will consider a compact Riemann surface $z_{0}(U)/\mathcal{M}(0)$
of genus $g$ with marking type $(h, g - h),$ induced by image
$z_{0}(\partial \mathcal{F}_{\tau_{0}}).$ In case $h \neq 0$ this
making we will add to homotopical basis, i.e.
$z_{0}(U)/\mathcal{M}(0) = [F_{0}, \{a_{k}, b_{k}\}^{g}_{k=1}],$
$a_{k}\cap b_{k} = O \in F_{0}.$ By the theorem 3 [22] we obtain
that $\tau_{0} = [F_{0}, \{\psi_{0}(a_{k}),
\psi_{0}(b_{k})\}^{g}_{k=1}],$ where $\psi_{0}$ is some
homeomorphism $F_{0}$ on itself, $\psi_{0}(O) = O.$

For continuous family of groups $\mathcal{M}(\xi) =
\{\widetilde{A}_{1,\xi},..., \widetilde{B}_{g,\xi}\}$ there exists
a continuous family  of quasiconformal mappings $f_{\xi}$ of plane
$\overline{\bf C}$ onto itself such that $\mathcal{M}(\xi) =
f_{\xi} \mathcal{M}(0) f^{-1}_{\xi},$ $\xi \in [0, 1].$ Hence,
there is a continuous family of quasiconformal homeomorphisms
$\widetilde{f_{\xi}}$ from surface $\tau_{0}$ on $\tau_{\xi} =
[\widetilde{f_{\xi}}(F_{0}),
\{\psi_{\xi}\widetilde{f_{\xi}}(a_{k}),
\psi_{\xi}\widetilde{f_{\xi}}(b_{k})\}^{g}_{k=1}],$ $\psi_{\xi} =
\widetilde{f_{\xi}}\psi_{0}(\widetilde{f_{\xi}})^{-1}$ is
homeomorphism $\widetilde{f_{\xi}}(F_{0})$ onto itself and
$\widetilde{f_{\xi}}$ map
$z_{0}(\mathcal{F}_{\tau_{0}})/\mathcal{M}[z_{0}]$ on
$f_{\xi}(z_{0}(\mathcal{F}_{\tau_{0}}))/\mathcal{M}(\xi).$ Let
$\tau_{\xi} = [F_{\xi},\widetilde{f_{\xi}}],$ where
$[F_{0},\widetilde{f_{0}}] \equiv [F_{0}, \{\psi_{0}(a_{k}),
\psi_{0}(b_{k})\}^{g}_{k=1}] = \tau_{0}$ and $\widetilde{f_{0}} =
1$ (identical mapping of $F_{0}).$ Notice that automorphism of
group $\pi_{1}(\widetilde{f_{\xi}}(F_{0}))$ with generators
$\{\widetilde{f_{\xi}}(a_{k}),
\widetilde{f_{\xi}}(b_{k})\}^{g}_{k=1},$ induced by $\psi_{\xi},$
coincides with automorphism of group $\pi_{1}(F_{0})$ with
generators $ \{a_{k},b_{k}\}^{g}_{k=1},$ induced by $\psi_{0}.$
Quasiconformal mapping $\widetilde{f_{\xi}}$ is lifted up to
quasiconformal mapping $\varphi_{\xi}: U \rightarrow U.$ So, as
$\widetilde{f_{\xi}}$ is determined up to  isotopy and
$\mathcal{M}(\xi)$ is continuous on $\xi,$ the mapping
$\varphi_{\xi}(t)$ is continuous on $U \times [0,1],$ and
appropriate normalized Fuchsian group $\Gamma_{\xi} =
\{A_{1,\xi},..., B_{g,\xi}\}$ has fundamental polygon
$\mathcal{F}_{\tau_{\xi}} =
\varphi_{\xi}(\mathcal{F}_{\tau_{0}}),$ $0 \leq \xi \leq 1.$
Moreover $\varphi_{0}$ is identical mapping of disc $U.$ We
received family of mappings
$$
z_{\xi} = f_{\xi} z_{0} \varphi_{\xi}^{-1}: U \rightarrow
z_{\xi}(U),0 \leq \xi \leq 1,
$$
such that  $z_{\xi}(\mathcal{F}_{\tau_{\xi}}) =
f_{\xi}(z_{0}(\mathcal{F}_{\tau_{0}})).$ Hence, we have a
continuous family of fundamental membranes, sides of whose are
paired identified by generators of the group $\mathcal{M}(\xi).$
Hence, $\{z_{\xi}\}$ is a continuous family locally schlicht
linearly polymorphic function on continuous family of marked
compact Riemann surfaces $[F_{\xi},
\{\psi_{\xi}\widetilde{f}_{\xi}(a_{k}),
\psi_{\xi}\widetilde{f}_{\xi}(b_{k})\}^{g}_{k=1}],$ $\xi \in [0,
1],$ respectively. Moreover,
$f_{\xi}(z_{0}(\mathcal{F}_{\tau_{0}}))$ defines a fundamental
membrane for locally schlicht analytic  linearly polymorphic
function  $z_{\xi}$ on $U$ with the monodromy group
$\mathcal{M}(\xi)$ [1]. Thus, we have constructed continuous path
$\widetilde{\gamma}_{0} = \{(\tau_{\xi}, q_{\xi} = \{z_{\xi},
t\}): \xi \in [0, 1]\}$ from $(\tau_{0}, q_{0})$ in $B_{\sigma}$
over $\{\mathcal{M}(\xi): \xi \in [0, 1]\}.$

Let now $\mathcal{M}[w] = \mathcal{M}(0)$ and $w = w(s)$
corresponds $([F_{w}, h_{w}], q_{w}) \in {\bf T}_{g} Q \backslash
 {\bf B}_{\sigma}$
with Fuchsian collection
$$
\{ F_{w},  \mathcal{F}_{w}, \Gamma_{w} = \{{A_{k,w}},
B_{k,w}\}_{k=1}^{g}\},
$$
where $\Gamma_{w}$ be Fuchsian group with fundamental polygon
$\mathcal{F}_{w},$ which uniformize $F_{w}$ in disc $U.$ By the
theorem of Kra-Gunning we have $w(U) = \overline{\bf C},$ i.e. $w$
is not uniformization for $[F_{w}, h_{w}].$ By the theorem 2 [1]
it is possible to achieve that appropriate vertices of membranes
$w(\mathcal{F}_{w})$ and $z_{0}(\mathcal{F}_{\tau_{0}})$ would
have identical coordinates. A fundamental membrane
$f_{\xi}(w(\mathcal{F}_{w}))$ defines at $\xi \in [0, 1]$ (not
analytic) locally schlicht linearly polymorphic function
$f_{\xi}w$ on $U$ such that $f_{\xi}w(A_{k,w}) =
\widetilde{A}_{k,\xi} f_{\xi} w,$ $f_{\xi} w (B_{k,w}) =
\widetilde{B}_{k,\xi} f_{\xi} w,$ $k = 1,..., g.$ Continuously
deforming the complex-analytic structure on these membranes with
the help quasiconformal mapping [1], [6, c.348], we receive that
$f_{\xi} w (\mathcal{F}_{w})$ defines a fundamental membrane for
analytic locally schlicht linearly polymorphic function on $U$
with the monodromy group $\mathcal{M}(\xi)$ for every $\xi \in [0,
1].$

Since $\gamma$ is continuously of $\xi,$ by the theorem 1[1],
using found the fundamental membranes for linearly polymorphic
function on $U,$ we  received a continuous path
$\widetilde{\gamma_{0}}$ from $(\tau_{0}, q_{0}) \in {\bf
B}_{\sigma}$ and $\widetilde{\gamma_{w}}$ from $([F_{w}, h_{w}],
q_{w}) \in
 {\bf T}_{g} Q \setminus {\bf B}_{\sigma}.$
Theorem is proved.

Remark 2.1. By the theorem 1[1] and the theorem 2.1 it follows
that:

1) mapping $p$ defined over space ${\bf T}(G_{\sigma}) = {\bf
T}(\Omega(G_{\sigma}), G_{\sigma}),$ as in the theorem, is a
topological covering;

2) since space ${\bf B}_{\sigma}$ is simply connected it is
universal covering space for ${\bf T}(G_{\sigma})$ and $p : {\bf
B}_{\sigma} \rightarrow {\bf T}(G_{\sigma})$ is the universal
covering mapping.

\vspace{3mm}

\centerline{{\bf 3. The exact variational formula for monodromy
group of the}}

\centerline{{\bf linear differential equation of the second order
and for the solution}}

\centerline{{\bf of the nonlinear Schwartz equation on compact
Riemann surface}}

\vspace{3mm}

In work [2] D.A. Hejhal have started the research of monodromy
group for linearly polymorphic function on compact Riemann surface
with the help of variational methods. He has found the first
variation for monodromy group. Then C. Earle [4] has deduced the
formula the first variation with the help quasiconformal mapping
of Riemann surfaces.

In this section we will received an exact variational formula for
monodromy group of the linear differential equation of the second
order and the first variation for solution of the Schwartz
equation on compact Riemann surface.

Let $F$ be a compact Riemann surface of genus $g, g \geq 2;
\pi_{1}(F, O)$ is the fundamental group for $F$ with basic point
$O,$ and $\Gamma$ is the group of covering transformations for
universal covering $(U, \pi)$ over $F.$ Here $U = \{t \in {\bf C}
: |t| < 1 \},$ $\pi : U \rightarrow F$ is natural projection and
$U/\Gamma = F.$ Fix a point $t_{0} \in  U,$ laying over $O,$ and
we will construct natural isomorphism of group $\pi_{1}(F, O)$ on
Fuchsian group $\Gamma$ of first kind.

A multivalued locally meromorphic function $z$ on $F,$  which
transforms linear fractionally under action group $\pi_{1}(F, O)$,
is called linearly polymorphic function on $F.$ Lifted it on $(U,
\pi),$ we received meromorphic single-valued function $z = z(t)$
on $U,$ which satisfies the condition
$$
z(L t) = \tilde{L} z(t), \tilde{L} \in PSL(2, {\bf C}), \eqno (1)
$$
for $L \in \Gamma, t \in U.$ Mapping $L \mapsto \tilde{L}$ gives a
homomorphism of group $\Gamma$ in group $PSL(2, {\bf C}).$ Group,
which consist of mappings $\tilde{L},$ when $L$ runs $\Gamma,$ is
called by monodromy group for function $z = z(t).$ Function
$$
2q(t) = \{z, t\} = (z''/z')' - \frac{1}{2}(z''/z')^{2}
$$
satisfies the relation
$$
q(t) = q(L t) L'(t)^{2}
$$
for $L \in \Gamma, t \in U.$ Hence, $q(t)$ defines the quadratic
differential on $F = U/\Gamma.$

A locally schlicht linearly polymorphic function  $z = z(t)$
satisfies on $U$ of the Schwartz equation
$$
\{z, t\} = 2q(t), \eqno (2)
$$
where $q(t)$ be holomorphic function on $U.$ Putting
$$
z(t) = v(t)/u(t), v(t) = z(t)/\sqrt{z'(t)}, u(t) = 1/\sqrt{z'(t)},
\eqno (3)
$$
we receive that $v = v(t), u = u(t)$ satisfy of the linear
differential equation of the second order
$$
u''(t) + q(t) u(t) = 0 \eqno (4)
$$
on $F = U/\Gamma.$ We denote by $Q(F)$ the vector space
holomorphic quadratic (abelian) differentials on $F = U/\Gamma.$
For equation (2) we will consider only the normalized solutions $z
= z(t)$ such that
$$
z(t) = (t - t_{0}) + O((t - t_{0})^{3}),
 \sqrt{z'(t)} = 1 + O((t - t_{0})^{2}), t \rightarrow t_{0}. \eqno (5)
$$
We obtain that $v = v(t)$ and $u = u(t)$ are two linearly
independent the solution of the equation (4) with conditions
$$
u(t_{0}) = 1, u'(t_{0}) = 0,
 v(t_{0}) = 0, v'(t_{0}) = 1, \eqno (6)
$$
for any $q = q(t) dt^{2} \in Q(F).$ It is well known that the
conditions (5) and (6) define the unique solutions of the equation
(2) and (4) respectively.

In work [2] was found the following relation
$$
\left(\begin{array}{l}
  {v(L t)} \\
  {u(L t)}
\end{array}\right)
 = \xi_{L}(t) \left(\begin{array}{ll}
  {\alpha} & {\beta} \\
  {\gamma} & {\delta}
\end{array}\right) \left(\begin{array}{l}
  {v(t)} \\
  {u(t)}
\end{array}\right), \eqno (7)
$$
where
$$
z(L t) = \tilde{L} z(t) = (\alpha z(t) + \beta) (\gamma z(t) +
\delta)^{-1}, \xi_{L}(t) = \sqrt{L'(t)},
$$
and for any $L, K \in \Gamma$ valid a relations $\xi_{L K}(t) =
\xi_{L}(K t) \xi_{K}(t).$ There is a special choise of a sign at
$\xi_{L}(t)$ and at a matrix $\left(\begin{array}{ll}
  {\alpha} & {\beta} \\
  {\gamma} & {\delta}
\end{array}\right)$
(see also [11, с. 163]).

Choose arbitrary $r = r(t) dt^{2}, q = q(t) dt^{2}$ from $Q(F).$
Consider normalized by (5) a solutions $z(t, h)$ of the Schwartz
equation
$$
\{z, t\} = 2[r(t) + h q(t)] \eqno (8)
$$
and normalized by (6) a solutions $v(t, h)$ and $u(t,h)$ of the
linear equation
$$
u''(t) + [r(t) + h q(t)] u(t) = 0, \eqno (9)
$$
where $h \in {\bf C}, |h| < \varepsilon,$ $ \varepsilon$ is
sufficiently small a positive number. Here $z(t, h) = v(t, h)/u(t,
h).$ Applying the Poincare theorem about small parameter [23], the
Caushy-Kovalevski theorem, we receive  Hartog's series [24]
$$
u(t, h) = u(t) + u_{1}(t) h + u_{2}(t) h^{2} + ... + u_{n}(t)
h^{n} + ... ,
$$
$$
v(t, h) = v(t) + v_{1}(t) h + v_{2}(t) h^{2} + ... , \eqno (10)
$$
which are uniformly converged on any compact in $U \times \{h :
|h| < \varepsilon \},$ where $u(t), v(t), u_{i}(t), v_{i}(t)$ are
holomorphic functions on $U.$ From normalization (6) it follows
that
$$
u(t_{0}) = 1, u'(t_{0}) = 0, v(t_{0}) = 0, v'(t_{0}) = 1,
$$
$$
u_{i}(t_{0}) = v_{i}(t_{0}) = u'_{i}(t_{0}) = v'_{i}(t_{0}) = 0, i
\geq 1. \eqno (11)
$$
Substituting the series (10) in (9), we receive a infinite system
of pairs of equations
$$
\{\begin{array}{l}
  u''(t) + r(t) u(t) = 0
 \\
 v''(t) + r(t) v(t) = 0
\end{array}
$$
$$
\{\begin{array}{l} u''_{1}(t) + r(t) u_{1}(t) = -q(t) u(t)
\\
v''_{1}(t) + r(t) v_{1}(t) = -q(t) v(t)
\end{array}
$$
$$...$$
$$
\{\begin{array}{l} u''_{n}(t) + r(t) u_{n}(t) = -q(t) u_{n-1}(t)
\\
v''_{n}(t) + r(t) v_{n}(t) = -q(t) v_{n-1}(t)
\end{array}
$$
$$...$$
Solving the Caushy problem with zero initial conditions in point
$t_{0}$ for system with pair $(u_{1}(t), v_{1}(t))$ on a method of
an elementary variation of parameters, we  find that
$$
\left(\begin{array}{l}
  {v_{1}(t)} \\
  {u_{1}(t)}
\end{array}\right)
 = \int_{t_{0}}^{t} \left(\begin{array}{ll}
  {-q u v} & {q v v} \\
  {-q u u} & {q u v}
\end{array}\right) d s \left(\begin{array}{l}
  {v(t)} \\
  {u(t)}
\end{array}\right),   \eqno (12)
$$
where $u = u(t), v = v(t)$ are the solutions the Caushy problem
with initial conditions $u(t_{0}) = v'(t_{0}) = 1, u'(t_{0}) =
v(t_{0}) = 0$ for first pairs of equations. By induction for any
$n, n \geq 1,$ we will receive a relations
$$
\left(\begin{array}{l}
  {v_{n}(t)} \\
  {u_{n}(t)}
\end{array}\right) = \int_{t_{0}}^{t} \left(\begin{array}{ll}
  {-q u v_{n-1}} & {q v v_{n-1}} \\
  {-q u u_{n-1}} & {q v u_{n-1}}
\end{array}\right) d s \left(\begin{array}{l}
  {v(t)} \\
  {u(t)}
\end{array}\right) \equiv
$$
$$
 \equiv A_{n-1}(t)\left(\begin{array}{l}
  {v(t)} \\
  {u(t)}
\end{array}\right) = \int_{t_{0}}^{t} A_{n-2}(s) A(s) d s \left(\begin{array}{l}
  {v(t)} \\
  {u(t)}
\end{array}\right),            \eqno  (13)
$$
where
$$
A(s) = q(s) \left(\begin{array}{ll}
  {-u(s) v(s)} & {v^{2}(s)} \\
  {-u^{2}(s)} & {u(s) v(s)}
\end{array}\right), A_{0}(t) = \int_{t_{0}}^{t} A(s) d s.
$$
It is easy to see that formula (13) gives two solutions of the
Caushy problem with zero initial conditions in point $t_{0}$ for
system with pair $(u_{n}(t), v_{n}(t)).$

Hence we receive the exact variational formula for solutions of
the linear equation (9)
$$
\left(\begin{array}{l}
  {v(t,h)} \\
  {u(t,h)}
\end{array}\right) = \left(\begin{array}{l}
  {v(t)} \\
  {u(t)}
\end{array}\right) + h \left(\begin{array}{l}
  {v_{1}(t)} \\
  {u_{1}(t)}
\end{array}\right) + h^{2} \left(\begin{array}{l}
  {v_{2}(t)} \\
  {u_{2}(t)}
\end{array}\right) + ... + h^{n} \left(\begin{array}{l}
  {v_{n}(t)} \\
  {u_{n}(t)}
\end{array}\right) + ... =
$$
$$
[\left(\begin{array}{ll}
  {1} & {0} \\
  {0} & {1}
\end{array}\right) +  h \int_{t_{0}}^{t} A(s) d s +
h^{2} \int_{t_{0}}^{t} A_{0}(s) A(s) d s + ... +
 h^{n} \int_{t_{0}}^{t} A_{n-2}(s) A(s) d s + ...] \cdot
 $$
 $$
\cdot \left(\begin{array}{l}
  {v(t)} \\
  {u(t)}
\end{array}\right)\equiv
\Omega _{t_{0}}^{t}h A(s) d s \cdot \left(\begin{array}{l}
  {v(t)} \\
  {u(t)}
\end{array}\right).
$$
The expression in square brackets is called the matrizant of $h
A(s).$

For finding of the exact variational formula of the monodromy
group for function $z(t, h)$ it is necessary deduce some
relations. By analogy with (7) we  write for any $h, \mid h \mid <
\varepsilon,$
$$
z(L t, h) = \tilde{L}_{h} z(t, h) = (\alpha_{L}(h) z(t, h) +
\beta_{L}(h))
 (\gamma_{L}(h) z(t, h) + \delta_{L}(h))^{-1},
$$
$$
\left(\begin{array}{l}
  {v(L t, h)} \\
  {u(L t, h)}
\end{array}\right) = \xi_{L}(t) \left(\begin{array}{ll}
  {\alpha_{L}(h)} & {\beta_{L}(h)} \\
  {\gamma_{L}(h)} & {\delta_{L}(h)}
\end{array}\right) \left(\begin{array}{l}
  {v(t, h)} \\
  {u(t, h)}
\end{array}\right). \eqno(14)
$$
By elementary transformations, using (7), we obtain the following
relations
$$
\int_{L t_{0}}^{L t} A(x) \left(\begin{array}{ll}
  {\alpha} & {\beta} \\
  {\gamma} & {\delta}
\end{array}\right) d x = \left(\begin{array}{ll}
  {\alpha} & {\beta} \\
  {\gamma} & {\delta}
\end{array}\right) \int_{t_{0}}^{t} A(s) d s,
$$
$$
\left(\begin{array}{ll}
  {\alpha} & {\beta} \\
  {\gamma} & {\delta}
\end{array}\right) = \left(\begin{array}{ll}
  {\alpha_{L}(0)} & {\beta_{L}(0)} \\
  {\gamma_{L}(0)} & {\delta_{L}(0)}
\end{array}\right). \eqno (15)
$$
From the formulas (12),(7),(15) we have the equality
$$
\left(\begin{array}{l}
  {v_{1}(L t)} \\
  {u_{1}(L t)}
\end{array}\right) =
$$
$$
 = \int_{t_{0}}^{L t} A(s) d s
\cdot \left(\begin{array}{l}
  {v(L t)} \\
  {u(L t)}
\end{array}\right)
 = \int_{t_{0}}^{L t} A(s) d s
\cdot \xi_{L}(t) \left(\begin{array}{ll}
  {\alpha} & {\beta} \\
  {\gamma} & {\delta}
\end{array}\right)
\left(\begin{array}{l}
  {v(t)} \\
  {u(t)}
\end{array}\right) =
$$
$$
 = \xi_{L}(t) \int_{t_{0}}^{L t_{0}} A(s) \left(\begin{array}{ll}
  {\alpha} & {\beta} \\
  {\gamma} & {\delta}
\end{array}\right) d s \cdot
\left(\begin{array}{l}
  {v(t)} \\
  {u(t)}
\end{array}\right) +
$$
$$
 + \xi_{L}(t)
\left(\begin{array}{ll}
  {\alpha} & {\beta} \\
  {\gamma} & {\delta}
\end{array}\right)
\int_{t_{0}}^{ t} A(s)
 d s \cdot
\left(\begin{array}{l}
  {v(t)} \\
  {u(t)}
\end{array}\right) =
$$
$$
 = \xi_{L}(t) \int_{t_{0}}^{L t_{0}} A(s) d s \cdot
 \left(\begin{array}{ll}
  {\alpha} & {\beta} \\
  {\gamma} & {\delta}
\end{array}\right)
\left(\begin{array}{l}
  {v(t)} \\
  {u(t)}
\end{array}\right) +
\xi_{L}(t) \left(\begin{array}{ll}
  {\alpha} & {\beta} \\
  {\gamma} & {\delta}
\end{array}\right)
\left(\begin{array}{l}
  {v_{1}(t)} \\
  {u_{1}(t)}
\end{array}\right). \eqno   (16)
$$
Next, the equality (13),(7),(15) give us
$$
 \left(\begin{array}{l}
  {v_{2}(L t)} \\
  {u_{2}(L t)}
\end{array}\right) = \xi_{L}(t) \int_{t_{0}}^{L t}
(\int_{t_{0}}^{x} A(s_{1}) d s_{1})A(x) \left(\begin{array}{ll}
  {\alpha} & {\beta} \\
  {\gamma} & {\delta}
\end{array}\right) d x \left(\begin{array}{l}
  {v(t)} \\
  {u(t)}
\end{array}\right) =
$$
$$
 = \xi_{L}(t) A_{1}(L t_{0}) \left(\begin{array}{ll}
  {\alpha} & {\beta} \\
  {\gamma} & {\delta}
\end{array}\right) \left(\begin{array}{l}
  {v(t)} \\
  {u(t)}
\end{array}\right) + \xi_{L}(t) A_{0}(L t_{0}) \left(\begin{array}{ll}
  {\alpha} & {\beta} \\
  {\gamma} & {\delta}
\end{array}\right)\left(\begin{array}{l}
  {v_{1}(t)} \\
  {u_{1}(t)}
\end{array}\right) + \eqno (17)
$$
$$
 + \xi_{L}(t) \left(\begin{array}{ll}
  {\alpha} & {\beta} \\
  {\gamma} & {\delta}
\end{array}\right) \left(\begin{array}{l}
  {v_{2}(t)} \\
  {u_{2}(t)}
\end{array}\right), $$
since
$$
\int_{L t_{0}}^{L t}(\int_{t_{0}}^{x}
 A(s_{1}) d s_{1}) A(x) \left(\begin{array}{ll}
  {\alpha} & {\beta} \\
  {\gamma} & {\delta}
\end{array}\right) d x =
$$
$$
 = \int_{t_{0}}^{ t}(\int_{t_{0}}^{L(s)}
 A(s_{1}) d s_{1}) A(L s)
 \left(\begin{array}{ll}
  {\alpha} & {\beta} \\
  {\gamma} & {\delta}
\end{array}\right)
d L(s) =
$$
$$
 = \int_{t_{0}}^{ t}(\int_{t_{0}}^{L(s)}
 A(s_{1}) d s_{1})
 \left(\begin{array}{ll}
  {\alpha} & {\beta} \\
  {\gamma} & {\delta}
\end{array}\right)
 A(s) d s =
$$
$$
 = (\int_{t_{0}}^{ L t_{0}}
 A(s_{1}) d s_{1})
 \left(\begin{array}{ll}
  {\alpha} & {\beta} \\
  {\gamma} & {\delta}
\end{array}\right)
\int_{t_{0}}^{t}A(s)d s +
$$
$$
 + \int_{t_{0}}^{ t}(\int_{t_{0}}^{s}
 A(L x_{1})
 \left(\begin{array}{ll}
  {\alpha} & {\beta} \\
  {\gamma} & {\delta}
\end{array}\right)
 d L (x_1))A(s)d s =
$$
$$
 = A_{0}(L t_{0}) \left(\begin{array}{ll}
  {\alpha} & {\beta} \\
  {\gamma} & {\delta}
\end{array}\right)A_{0}(t) + \left(\begin{array}{ll}
  {\alpha} & {\beta} \\
  {\gamma} & {\delta}
\end{array}\right) A_{1}(t).
$$
By induction on $n, n \geq 2,$ we find the relations
$$
\left(\begin{array}{l}
  {v_{n}(L t)} \\
  {u_{n}(L t)}
\end{array}\right) = \int_{t_{0}}^{L t} q(s) \left(\begin{array}{ll}
  {v_{n-1}(s)} & {0} \\
  {u_{n-1}(s)} & {0}
\end{array}\right) \left(\begin{array}{ll}
  {-u(s)} & {v(s)} \\
  {0}     &  {0}
\end{array}\right) d s \left(\begin{array}{l}
  {v(L t)} \\
  {u(L t)}
\end{array}\right) =
$$
$$
 = \xi_{L}(t) A_{n-1}(L t_{0})
\left(\begin{array}{ll}
  {\alpha} & {\beta} \\
  {\gamma} & {\delta}
\end{array}\right)
 \left(\begin{array}{l}
  {v(t)} \\
  {u(t)}
\end{array}\right) +
$$
$$
 + \xi_{L}(t) \int_{L t_{0}}^{L t}
q(x)
 \left(\begin{array}{ll}
  {v_{n-1}(x)} & {0} \\
  {u_{n-1}(x)} & {0}
\end{array}\right)
\left(\begin{array}{ll}
  {-u(x)} & {v(x)} \\
  {0}     &   {0}
\end{array}\right)
\left(\begin{array}{ll}
  {\alpha} & {\beta} \\
  {\gamma} & {\delta}
\end{array}\right)
 d x
 \left(\begin{array}{l}
  {v(t)} \\
  {u(t)}
\end{array}\right) =
$$
$$
 = [\xi_{L}(t) A_{n-1}(L t_{0}) \left(\begin{array}{ll}
  {\alpha} & {\beta} \\
  {\gamma} & {\delta}
\end{array}\right) +
$$
$$
 + \xi_{L}(t) \int_{t_{0}}^{t}
\frac{q(s) \xi_{L}(s)} {L'(s)} \left(\begin{array}{ll}
  {v_{n-1}(L s)} & {0} \\
  {u_{n-1}(L s)} & {0}
\end{array}\right) \left(\begin{array}{ll}
  {-u(s)} & {v(s)} \\
  {0}     &   {0}
\end{array}\right) d s] \left(\begin{array}{l}
  {v(t)} \\
  {u(t)}
\end{array}\right) =
$$
$$
 = \xi_{L}(t) A_{n-1}(L t_{0}) \left(\begin{array}{ll}
  {\alpha} & {\beta} \\
  {\gamma} & {\delta}
\end{array}\right) \left(\begin{array}{l}
  {v(t)} \\
  {u(t)}
\end{array}\right) +
$$
$$
 + \xi_{L}(t) \int_{t_{0}}^{t} \frac{q(s) \xi_{L}^{2}(s)} {L'(s)}
[A_{n-2}(L t_{0}) \left(\begin{array}{ll}
  {\alpha} & {\beta} \\
  {\gamma} & {\delta}
\end{array}\right) \left(\begin{array}{ll}
  {v(s)} & {0} \\
  {u(s)} & {0}
\end{array}\right) +
$$
$$
 + \sum_{j=1}^{n-2} A_{n-2-j}(L t_{0}) \left(\begin{array}{ll}
  {\alpha} & {\beta} \\
  {\gamma} & {\delta}
\end{array}\right) \left(\begin{array}{ll}
  {v_{j}(s)} & {0} \\
  {u_{j}(s)} & {0}
\end{array}\right) +
$$
$$
+ \left(\begin{array}{ll}
  {\alpha} & {\beta} \\
  {\gamma} & {\delta}
\end{array}\right) \left(\begin{array}{ll}
  {v_{n-1}(s)} & {0} \\
  {u_{n-1}(s)} & {0}
\end{array}\right)] \left(\begin{array}{ll}
  {-u(s)} & {v(s)} \\
  {0}     &   {0}
\end{array}\right) d s \left(\begin{array}{l}
  {v(t)} \\
  {u(t)}
\end{array}\right) =
$$
$$
 = \xi_{L}(t) [A_{n-1}(L t_{0}) \left(\begin{array}{ll}
  {\alpha} & {\beta} \\
  {\gamma} & {\delta}
\end{array}\right) \left(\begin{array}{l}
  {v(t)} \\
  {u(t)}
\end{array}\right) +
$$
$$
 + \sum_{j=1}^{n-1} A_{n-1-j}(L t_{0}) \left(\begin{array}{ll}
  {\alpha} & {\beta} \\
  {\gamma} & {\delta}
\end{array}\right) \left(\begin{array}{l}
  {v_{j}(t)} \\
  {u_{j}(t)}
\end{array}\right) + \left(\begin{array}{ll}
  {\alpha} & {\beta} \\
  {\gamma} & {\delta}
\end{array}\right) \left(\begin{array}{l}
  {v_{n}(t)} \\
  {u_{n}(t)}
\end{array}\right)].       \eqno (18)
$$
Using the relations (14) and (16) - (18), we  receive
$$
\xi_{L}(t) \left(\begin{array}{ll}
  {\alpha_{L}(h)} & {\beta_{L}(h)} \\
  {\gamma_{L}(h)} & {\delta_{L}(h)}
\end{array}\right)
\left(\begin{array}{l}
  {v(t, h)} \\
  {u(t, h)}
\end{array}\right) =
\left(\begin{array}{l}
  {v(L t, h)} \\
  {u(L t, h)}
\end{array}\right) =
$$
$$
 = \left(\begin{array}{l}
  {v(L t)} \\
  {u(L t)}
\end{array}\right) + h \left(\begin{array}{l}
  {v_{1}(L t)} \\
  {u_{1}(L t)}
\end{array}\right) + h^{2} \left(\begin{array}{l}
  {v_{2}(L t)} \\
  {u_{2}(L t)}
\end{array}\right) + ... +
 h^{n} \left(\begin{array}{l}
  {v_{n}(L t)} \\
  {u_{n}(L t)}
\end{array}\right) + ... =
$$
$$
 = \xi_{L}(t) [\left(\begin{array}{ll}
  {\alpha} & {\beta} \\
  {\gamma} & {\delta}
\end{array}\right) \left(\begin{array}{l}
  {v(t)} \\
  {u(t)}
\end{array}\right) + h A_{0}(L t_{0}) \left(\begin{array}{ll}
  {\alpha} & {\beta} \\
  {\gamma} & {\delta}
\end{array}\right) \left(\begin{array}{l}
  {v(t)} \\
  {u(t)}
\end{array}\right) +
$$
$$
 + h \left(\begin{array}{ll}
  {\alpha} & {\beta} \\
  {\gamma} & {\delta}
\end{array}\right) \left(\begin{array}{l}
  {v_{1}(t)} \\
  {u_{1}(t)}
\end{array}\right) + h^{2} A_{1}(L t_{0}) \left(\begin{array}{ll}
  {\alpha} & {\beta} \\
  {\gamma} & {\delta}
\end{array}\right) \left(\begin{array}{l}
  {v(t)} \\
  {u(t)}
\end{array}\right) +
$$
$$
 + h^{2} A_{0}(L t_{0}) \left(\begin{array}{ll}
  {\alpha} & {\beta} \\
  {\gamma} & {\delta}
\end{array}\right) \left(\begin{array}{l}
  {v_{1}(t)} \\
  {u_{1}(t)}
\end{array}\right) + h^{2} \left(\begin{array}{ll}
  {\alpha} & {\beta} \\
  {\gamma} & {\delta}
\end{array}\right) \left(\begin{array}{l}
  {v_{2}(t)} \\
  {u_{2}(t)}
\end{array}\right) + ... +
$$
$$
 + h^{n} A_{n-1}(L t_{0}) \left(\begin{array}{ll}
  {\alpha} & {\beta} \\
  {\gamma} & {\delta}
\end{array}\right) \left(\begin{array}{l}
  {v(t)} \\
  {u(t)}
\end{array}\right) + h^{n} \sum_{j=1}^{n-1} A_{n-1-j}(L t_{0}) \left(\begin{array}{ll}
  {\alpha} & {\beta} \\
  {\gamma} & {\delta}
\end{array}\right) \left(\begin{array}{l}
  {v_{j}(t)} \\
  {u_{j}(t)}
\end{array}\right) +
$$
$$
 + h^{n} \left(\begin{array}{ll}
  {\alpha} & {\beta} \\
  {\gamma} & {\delta}
\end{array}\right) \left(\begin{array}{l}
  {v_{n}(t)} \\
  {u_{n}(t)}
\end{array}\right) + ...] =
$$
$$
 = \xi_{L}(t) [\left(\begin{array}{ll}
  {\alpha} & {\beta} \\
  {\gamma} & {\delta}
\end{array}\right) \left(\begin{array}{l}
  {v(t, h)} \\
  {u(t, h)}
\end{array}\right) + h A_{0}(L t_{0}) \left(\begin{array}{ll}
  {\alpha} & {\beta} \\
  {\gamma} & {\delta}
\end{array}\right) \left(\begin{array}{l}
  {v(t, h)} \\
  {u(t, h)}
\end{array}\right) +
$$
$$
 + h^{2} A_{1}(L t_{0}) \left(\begin{array}{ll}
  {\alpha} & {\beta} \\
  {\gamma} & {\delta}
\end{array}\right) \left(\begin{array}{l}
  {v(t, h)} \\
  {u(t, h)}
\end{array}\right) + ... +
 h^{n} A_{n-1}(L t_{0})
 \left(\begin{array}{ll}
  {\alpha} & {\beta} \\
  {\gamma} & {\delta}
\end{array}\right)
\left(\begin{array}{l}
  {v(t, h)} \\
  {u(t, h)}
\end{array}\right) + ...].
$$
Hence, we have the exact variational formula for elements of
monodromy group of function $z(t, h)$
$$
\left(\begin{array}{ll}
  {\alpha_{L}(h)} & {\beta_{L}(h)} \\
  {\gamma_{L}(h)} & {\delta_{L}(h)}
\end{array}\right) =
$$
$$
 = [\left(\begin{array}{ll}
  {1} & {0} \\
  {0} & {1}
\end{array}\right)
 + h A_{0}(L t_{0}) + h^{2} A_{1}(L t_{0}) + ... +
 h^{n} A_{n-1}(L t_{0}) + ...]
 \left(\begin{array}{ll}
  {\alpha_{L}(0)} & {\beta_{L}(0)} \\
  {\gamma_{L}(0)} & {\delta_{L}(0)}
\end{array}\right)
 \equiv
$$
$$
 [\Omega_{t_{0}}^{L t_{0}}h A(s) d s]
\left(\begin{array}{ll}
  {\alpha_{L}(0)} & {\beta_{L}(0)} \\
  {\gamma_{L}(0)} & {\delta_{L}(0)}
\end{array}\right),
$$
where
$$
 A_{0}(x) = \int_{t_{0}}^{x} A(s) d s, A(s) = q(s)
 \left(\begin{array}{ll}
  {-u(s) v(s)} & {v^{2}(s)} \\
  {-u^{2}(s)} & {u(s) v(s)}
\end{array}\right),
$$
$$
A_{n}(x) = \int_{t_{0}}^{x} A_{n-1}(s) A(s) d s, n \geq 1, u(s) =
u(s, 0), v(s) = v(s, 0).
$$
Now we will deduce the variational formula for the solution of the
Schwartz equation (8). From the formula (10) we have
$$
z(t, h) = \frac{v(t, h)}{u(t, h)} = \frac{v(t) + v_{1}(t) h +
v_{2}(t) h^{2} + ...} {u(t) + u_{1}(t) h + u_{2}(t) h^{2} + ...} =
$$
$$
 = \frac{v(t)}{u(t)} + h [\frac{v_{1}(t) u(t) - v(t) u_{1}(t)}
{u^{2}(t)}] + o(h), h \rightarrow 0.
$$
Using (12), we  receive the equality
$$
\frac{v_{1}(t) u(t) - v(t) u_{1}(t)}{u^{2}(t)} =
$$
$$
 = \frac{1}{u^{2}(t)}[u(t)(v(t) \int_{t_{0}}^{t}
(-q u v) d s + u(t) \int_{t_{0}}^{t} q v^{2} d s) - v(t)(v(t)
\int_{t_{0}}^{t} (-q u^{2}) d s +
$$
$$
 + u(t) \int_{t_{0}}^{t} q v u d s)] =
$$
$$
 = \frac{1}{u^{2}(t)}[u^{2}(t) \int_{t_{0}}^{t} q v^{2} d s -
 2 u(t) v(t) \int_{t_{0}}^{t} q u v d s +
v^{2}(t) \int_{t_{0}}^{t} q u^{2} d s] =
$$
$$
 = \int_{t_{0}}^{t} q(s) [v(s) - z(t, 0) u(s)]^{2} d s,
$$
where $u = u(s, 0), v = v(s, 0).$ Hence, we have the variational
formula
$$
z(t, h) = z(t, 0) + h \int_{t_{0}}^{t} q(s)[v(s) - z(t, 0)
u(s)]^{2} d s + o(h), h \rightarrow 0.
$$
By applying the standard formulas for coefficients  can be
received any variational term for $z(t, h).$

Consider the Schwartz equation
$$
\{z, t \} = 2[r(t) + \sum_{j=1}^{d} h_{j} q_{j}(t)] \eqno(19)
$$
and the linear equation
$$
u''(t) + [r(t) + \sum_{j=1}^{d} h_{j} q_{j}(t)] u(t) = 0 \eqno(20)
$$
on $F = U/\Gamma,$ where $q_{1}(t) d t^{2},..., q_{d}(t) d t^{2}$
is a basis in space $Q(F),$ $h = (h_{1},..., h_{d}) \in {\bf
C}^{d},$ $d = 3g - 3.$ Again we will consider only normalized
solutions $z(t, h)$ and $v(t, h),$ $u(t, h),$ with conditions (5)
and (6) respectively, for  any $h$ such that $|h| = \max_{\{1 \leq
j \leq d \}} |h_{j}| < \varepsilon, \varepsilon$ is sufficiently
small a positive number.

Now we deduce the variational formulas for the differential
equations (19) and (20). By the Poincare theorem about small
parameter [23], using decomposition for power series on
homogeneous polynomials in polydisk $ \{ t \in {\bf C} : |t| <
\delta\} \times \{h \in {\bf C}^{3g-3} : |h| < \varepsilon\}$ [24,
p. 52] and the analytical continuation, we will receive the series
$$
u(t, h) = \sum_{|k|=0}^{\infty} u_{|k|;(k_{1},..., k_{d})}(t)
h_{1}^{k_{1}}...h_{d}^{k_{d}} \equiv \sum_{|k|=0}^{\infty}
u_{|k|;k}(t) h^k,
$$
$$
v(t, h) =
 \sum_{|k|=0}^{\infty} v_{|k|;k}(t)
h^k = v(t) + \sum_{|k|=1}v_{1;k}(t)h^k +
$$
$$
\sum_{|k|=2}v_{2;k}(t)h^k + ... + \sum_{|k|=n}v_{n;k}(t)h^k + ...,
\eqno (21)
$$
which are uniformly converged on any compact in polydisk $ \{ t :
|t| < 1 \} \times \{h : |h| < \varepsilon \}.$ Here $k =
(k_{1},..., k_{d})$ is the vector with integer nonnegative
coordinates, and $|k| = k_{1} + ... + k_{d}.$ From normalization
$$
u(t_{0}, h) = 1 = v'(t_{0}, h), u'(t_{0}, h) = 0 = v(t_{0}, h)
$$
for any $h,$ we  receive
$$
u_{n;k}(t_{0}) =
u_{n;k}'(t_{0}) = v_{n;k}(t_{0}) = v_{n;k}'(t_{0}) = 0, n \geq 1,
|k| = n. \eqno (22)
$$
Substituting the series (21) in the equation (20), we receive a
infinite system of pairs of linear differential equations
$$
\{\begin{array}{l}
  u''(t) + r(t) u(t) = 0
 \\
 v''(t) + r(t) v(t) = 0,
\end{array}
$$
$$
\{\begin{array}{l} u''_{1;k}(t) + r(t) u_{1;k}(t) = -q_{j}(t) u(t)
\\
v''_{1;k}(t) + r(t) v_{1;k}(t) = -q_{j}(t) v(t),
\end{array}
$$
for $k = e_{j} = (0,..., 0, 1, 0,..., 0),$ where 1 stands on the
$j$-th place, $j = 1,..., d;$
$$...$$
$$
\{\begin{array}{l} u''_{n;k}(t) + r(t) u_{n;k}(t) = -\sum_{\{j:
k_{j}\neq 0\}}q_{j}(t) u_{n-1;k-e_{j}}(t)
\\
v''_{n;k}(t) + r(t) v_{n;k}(t) = -\sum_{\{j: k_{j}\neq
0\}}q_{j}(t) v_{n-1;k-e_{j}}(t),
\end{array}
$$
for $k = (k_{1},..., k_{d}), |k| = n;$
$$
...
$$
Hence, by method of elementary variation of parameters we  find
$$
\left(\begin{array}{l}
  {v_{1;k}(t)} \\
  {u_{1;k}(t)}
\end{array}\right)
 = \int_{t_{0}}^{t} \left(\begin{array}{ll}
  {-q_{j} u v} & {q_{j} v v} \\
  {-q_{j} u u} & {q_{j} u v}
\end{array}\right) d s \left(\begin{array}{l}
  {v(t)} \\
  {u(t)}
\end{array}\right)
 =
$$
$$
 = \int_{t_{0}}^{t}
M_{j}(s)d s \left(\begin{array}{l}
  {v(t)} \\
  {u(t)}
\end{array}\right)
\equiv A_{0;k}(t) \left(\begin{array}{l}
  {v(t)} \\
  {u(t)}
\end{array}\right),
$$
for $k = e_{j}$ and $M_{j} = q_{j}\left(\begin{array}{ll}
  -u v &  v v \\
  -u u &  u v
\end{array}\right),$
$j = 1,..., d.$

Then for $k = (k_{1},..., k_{d}), |k| = 2,$ we have
$$
\left(\begin{array}{l}
  {v_{2;k}(t)} \\
  {u_{2;k}(t)}
\end{array}\right) = \sum_{\{j:k_{j}\neq 0 \}} \int_{t_{0}}^{t} \left(\begin{array}{ll}
  {-q_{j} u v_{1;k-e_{j}}} & {q_{j} v v_{1;k-e_{j}}} \\
  {-q_{j} u u_{1;k-e_{j}}} & {q_{j} v u_{1;k-e_{j}}}
\end{array}\right) d s \left(\begin{array}{l}
  {v(t)} \\
  {u(t)}
\end{array}\right)
 =
 $$
 $$
 = \sum_{\{j:k_{j}\neq 0\}}
\int_{t_{0}}^{t}A_{0;k-e_{j}}(s) M_{j}(s)d s
\left(\begin{array}{l}
  {v(t)} \\
  {u(t)}
\end{array}\right)
\equiv
 A_{1;k}(t)\left(\begin{array}{l}
  {v(t)} \\
  {u(t)}
\end{array}\right).
$$
Similarity for $k = (k_{1},..., k_{d}), |k| = n,$ we will receive
equality
$$
\left(\begin{array}{l}
  {v_{n;k}(t)} \\
  {u_{n;k}(t)}
\end{array}\right) = \sum_{\{j:k_{j}\neq 0 \}} \int_{t_{0}}^{t} \left(\begin{array}{ll}
  {-q_{j} u v_{n-1;k-e_{j}}} & {q_{j} v v_{n-1;k-e_{j}}} \\
  {-q_{j} u u_{n-1;k-e_{j}}} & {q_{j} v u_{n-1;k-e_{j}}}
\end{array}\right) d s \left(\begin{array}{l}
  {v(t)} \\
  {u(t)}
\end{array}\right)
 =
 $$
 $$
 = \sum_{\{j:k_{j}\neq 0 \}} \int_{t_{0}}^{t}A_{n-2;k-e_{j}}(s) M_{j}(s)d s
\left(\begin{array}{l}
  {v(t)} \\
  {u(t)}
\end{array}\right)
\equiv A_{n-1;k}(t)\left(\begin{array}{l}
  {v(t)} \\
  {u(t)}
\end{array}\right).
$$

Thus, we have proved the following

{\bf Theorem 3.1.} \textsl{For pair of independent solutions of
the equation (20) with normalization (6) valid the exact
variational formula
$$
\left(\begin{array}{l}
  {v(t,h)} \\
  {u(t,h)}
\end{array}\right) =
\left(\begin{array}{l}
  {v(t)} \\
  {u(t)}
\end{array}\right) +
 \sum_{|k|=1}\left(\begin{array}{l}
  {v_{1;k}(t)} \\
  {u_{1;k}(t)}
\end{array}\right)
h^{k} +
$$
$$
 + \sum_{|k|=2}\left(\begin{array}{l}
  {v_{2;k}(t)} \\
  {u_{2;k}(t)}
\end{array}\right)
h^{k} + ... + \sum_{|k|=n}\left(\begin{array}{l}
  {v_{n;k}(t)} \\
  {u_{n;k}(t)}
\end{array}\right)
h^{k} + ... =
$$
$$
 [\left(\begin{array}{ll}
  {1} & {0} \\
  {0} & {1}
\end{array}\right) +
\sum_{|k|=1}A_{0;k}(t)h^{k} +
 \sum_{|k|=2}A_{1;k}(t)h^{k} +
... +  \sum_{|k|=n}A_{n-1;k}(t)h^{k} + ...] \left(\begin{array}{l}
  {v(t)} \\
  {u(t)}
\end{array}\right),
$$
where $A_{0;k}(t) = \int_{t_{0}}^{t}M_{j}(s) ds$ for $k = e_{j},$
$M_{j}(s) = q_{j}(s) \left(\begin{array}{ll}
  {-u(s)v(s)} & {v(s)v(s)} \\
  {-u(s)u(s)} & {u(s)v(s)}
\end{array}\right),$
$j = 1,..., d;$
$$
A_{n;(k_{1},..., k_{d})}(t) = \int_{t_{0}}^{t} \sum_{\{j:k_{j}\neq
0\}} A_{n-1;k-e_{j}}(s)M_{j}(s) ds,
$$
$n \geq 1, |k| \geq 2, d = 3g - 3, |h| < \varepsilon.$}

Derive the first variation for the solution of the Schwartz
equation (19):
$$
z(t, h) = \frac{v(t, h)}{u(t, h)} = \frac{v(t) + \sum _{i=1}^{d}
v_{i}(t) h_{i} + ...} {u(t) + \sum_{i=1}^{d} u_{i}(t) h_{i} + ...}
=
$$
$$
 = (v(t) + \sum_{i=1}^{d} v_{i}(t) h_{i} + ...)
(\frac{1}{u(t)} + \sum_{i=1}^{d} h_{i}
(-\frac{u_{i}(t)}{u^{2}(t)}) + ...) =
$$
$$
 = z(t, 0) + \sum_{i=1}^{d} h_{i}[\frac{u(t) v_{i}(t) - v(t)u_{i}(t)}{u^{2}(t)}] + ... =
$$
$$
 = z(t, 0) + \sum_{i=1}^{d} h_{i} \int_{t_{0}}^{t} q_{i}(s)(v(s) - z(t, 0) u(s))^{2} ds +
...
$$

Thus we have proved the following

{\bf Theorem 3.2.} {\sl{Let $z(t, h)$ be the solution of the
Schwartz equation (19) with normalization (5). Then  for $h =
(h_{1},..., h_{d}), |h| < \varepsilon,$ valid the variational
formula
$$
z(t, h) = z(t, 0) + \sum_{i=1}^{d} h_{i} \int_{t_{0}}^{t}
q_{i}(s)(v(s) - z(t, 0) u(s))^{2} ds +
 o(|h|), |h| \rightarrow 0.
$$}}

Remark 3.1. Applying the formulas for coefficients quotient power
series can be receive also any variational term for the solution
$z(t, h)$ of the Schwartz equation (19) on fixed compact Riemann
surface $F = U/\Gamma.$

To conclude of the variational formulas for elements of monodromy
groups it is necessary to find some additional relations connected
with group $\Gamma.$ First of all, for $k = (k_{1},..., k_{d}) =
e_{j}, |k| = 1,$ $j = 1,..., d,$ we find
$$
\left(\begin{array}{l}
  {v_{1;k}(L t)} \\
  {u_{1;k}(L t)}
\end{array}\right) =
\int_{t_{0}}^{L t} M_{j}(s) d s
 \left(\begin{array}{l}
  {v(L t)} \\
  {u(L t)}
\end{array}\right)
 =
 $$
 $$
 = \xi_{L}(t) \int_{t_{0}}^{L t_{0}} M_{j}(s) \left(\begin{array}{ll}
  {\alpha} & {\beta} \\
  {\gamma} & {\delta}
\end{array}\right) d s
\left(\begin{array}{l}
  {v(t)} \\
  {u(t)}
\end{array}\right) +
$$
$$
 + \xi_{L}(t) \int_{L t_{0}}^{L t} M_{j}(s) \left(\begin{array}{ll}
  {\alpha} & {\beta} \\
  {\gamma} & {\delta}
\end{array}\right) d s
\left(\begin{array}{l}
  {v(t)} \\
  {u(t)}
\end{array}\right)
 =
\xi_{L}(t) A_{0;k} (L t_{0}) \left(\begin{array}{ll}
  {\alpha} & {\beta} \\
  {\gamma} & {\delta}
\end{array}\right)
\left(\begin{array}{l}
  {v(t)} \\
  {u(t)}
\end{array}\right)
$$
$$
 + \xi_{L}(t)
\left(\begin{array}{ll}
  {\alpha} & {\beta} \\
  {\gamma} & {\delta}
\end{array}\right)
\int_{t_{0}}^{ t} M_{j}(s)
 d s
\left(\begin{array}{l}
  {v(t)} \\
  {u(t)}
\end{array}\right) =
$$
$$
 = \xi_{L}(t) A_{0;k}(L t_{0})
 \left(\begin{array}{ll}
  {\alpha} & {\beta} \\
  {\gamma} & {\delta}
\end{array}\right)
\left(\begin{array}{l}
  {v(t)} \\
  {u(t)}
\end{array}\right) +
\xi_{L}(t) \left(\begin{array}{ll}
  {\alpha} & {\beta} \\
  {\gamma} & {\delta}
\end{array}\right)
\left(\begin{array}{l}
  {v_{1;k}(t)} \\
  {u_{1;k}(t)}
\end{array}\right).
$$
For $k = (k_{1},..., k_{d}), |k| = 2,$ we get
$$
\left(\begin{array}{l}
  {v_{2;k}(L t)} \\
  {u_{2;k}(L t)}
\end{array}\right) = \xi_{L}(t) A_{1;k}(L t)
\left(\begin{array}{ll}
  {\alpha} & {\beta} \\
  {\gamma} & {\delta}
\end{array}\right) \left(\begin{array}{l}
  {v(t)} \\
  {u(t)}
\end{array}\right) =
$$

$$
 = \xi_{L}(t) A_{1;k}(L t_{0}) \left(\begin{array}{ll}
  {\alpha} & {\beta} \\
  {\gamma} & {\delta}
\end{array}\right)
\left(\begin{array}{l}
  {v(t)} \\
  {u(t)}
\end{array}\right)
 +
 $$
 $$
 + \xi_{L}(t) \int_{L (t_{0})}^{L(t)} \sum _{\{j:k_{j}\neq 0\}}
A_{0;k-e_{j}}(s) M_{j}(s) ds
 \left(\begin{array}{ll}
  {\alpha} & {\beta} \\
  {\gamma} & {\delta}
\end{array}\right)\left(\begin{array}{l}
  {v(t)} \\
  {u(t)}
\end{array}\right) =
$$
$$
 = \xi_{L}(t) A_{1;k}(L t_{0}) \left(\begin{array}{ll}
  {\alpha} & {\beta} \\
  {\gamma} & {\delta}
\end{array}\right) \left(\begin{array}{l}
  {v(t)} \\
  {u(t)}
\end{array}\right) +
$$
$$
 + \xi_{L}(t)\sum_{\{j:k_{j}\neq 0\}}
 A_{0;k-e_{j}}(L t_{0})
\left(\begin{array}{ll}
  {\alpha} & {\beta} \\
  {\gamma} & {\delta}
\end{array}\right) A_{0;e_{j}}(t)
\left(\begin{array}{l}
  {v(t)} \\
  {u(t)}
\end{array}\right) +
$$
$$
 +
 \xi_{L}(t) \left(\begin{array}{ll}
  {\alpha} & {\beta} \\
  {\gamma} & {\delta}
\end{array}\right)
\sum_{\{j:k_{j}\neq 0\}} \int_{t_{0}}^{t}A_{0;k-e_{j}}(s)
M_{j}(s)ds \left(\begin{array}{l}
  {v(t)} \\
  {u(t)}
\end{array}\right) =
$$
$$
 = \xi_{L}(t) A_{1;k}(L t_{0})
 \left(\begin{array}{ll}
  {\alpha} & {\beta} \\
  {\gamma} & {\delta}
\end{array}\right)
\left(\begin{array}{l}
  {v(t)} \\
  {u(t)}
\end{array}\right) +
$$
$$
 + \xi_{L}(t)\sum_{\{j:k_{j}\neq 0 \}}
 A_{0;k-e_{j}}(L t_{0})
\left(\begin{array}{ll}
  {\alpha} & {\beta} \\
  {\gamma} & {\delta}
\end{array}\right)
\left(\begin{array}{l}
  {v_{1;e_{j}}(t)} \\
  {u_{1;e_{j}}(t)}
\end{array}\right) +
$$
$$
 + \xi_{L}(t)
\left(\begin{array}{ll}
  {\alpha} & {\beta} \\
  {\gamma} & {\delta}
\end{array}\right)
\left(\begin{array}{l}
  {v_{2;k}(t)} \\
  {u_{2;k}(t)}
\end{array}\right).
$$
For $k = ( k_{1} ,..., k_{d} ),$ $ |k| = 3,$ we  find
$$
\left(\begin{array}{l}
  {v_{3;k}(L t)} \\
  {u_{3;k}(L t)}
\end{array}\right) =
 \xi_{L}(t) A_{2;k}(L t_{0})
 \left(\begin{array}{ll}
  {\alpha} & {\beta} \\
  {\gamma} & {\delta}
\end{array}\right)
\left(\begin{array}{l}
  {v(t)} \\
  {u(t)}
\end{array}\right) +
$$
$$
 + \xi_{L}(t) \int_{L(t_{0})}^{L(t)} \sum _{\{j:k_{j}\neq 0\}}
A_{1;k-e_{j}}(s) M_{j}(s) ds \left(\begin{array}{ll}
  {\alpha} & {\beta} \\
  {\gamma} & {\delta}
\end{array}\right)
\left(\begin{array}{l}
  {v(t)} \\
  {u(t)}
\end{array}\right) =
$$
$$
 = \xi_{L}(t) A_{2;k}(L t_{0})
\left(\begin{array}{ll}
  {\alpha} & {\beta} \\
  {\gamma} & {\delta}
\end{array}\right)
\left(\begin{array}{l}
  {v(t)} \\
  {u(t)}
\end{array}\right) +
$$
$$
 + \xi_{L}(t)\sum_{\{j:k_{j}\neq 0\}}
 A_{1;k-e_{j}}(L t_{0})
\left(\begin{array}{ll}
  {\alpha} & {\beta} \\
  {\gamma} & {\delta}
\end{array}\right)
\left(\begin{array}{l}
  {v_{1;e_{j}}(t)} \\
  {u_{1;e_{j}}(t)}
\end{array}\right) +
$$
$$
 +
 \xi_{L}(t) \sum_{\{j,i:k_{j}\neq 0, (k-e_{j})_{i}\neq 0\}}
A_{0;k-e_{j}-e_{i}}(L t_{0}) \left(\begin{array}{ll}
  {\alpha} & {\beta} \\
  {\gamma} & {\delta}
\end{array}\right)
\left(\begin{array}{l}
  {v_{2;e_{j}+e_{i}}(t)} \\
  {u_{2;e_{j}+e_{i}}(t)}
\end{array}\right) +
$$
$$
+ \xi_{L}(t) \left(\begin{array}{ll}
  {\alpha} & {\beta} \\
  {\gamma} & {\delta}
\end{array}\right)
\left(\begin{array}{l}
  {v_{3;k}(t)} \\
  {u_{3;k}(t)}
\end{array}\right).
$$

By induction for any $n > 1,$ we  receive the equality
$$
\left(\begin{array}{l}
  {v_{n + 1;k}(L t)} \\
  {u_{n + 1;k}(L t)}
\end{array}\right) =
$$
$$
 = \sum_{\{j:k_{j}\neq 0\}}\int_{t_{0}}^{L t} q_{j}(s)
\left(\begin{array}{ll}
  {v_{n;k-e_{j}}(s)} & {0} \\
  {u_{n;k-e_{j}}(s)} & {0}
\end{array}\right)
\left(\begin{array}{ll}
  {-u(s)} & {v(s)} \\
  {0} & {0}
\end{array}\right)
 ds
\left(\begin{array}{l}
  {v(L t)} \\
  {u(L t)}
\end{array}\right) =
$$
$$
 = \xi_{L}(t) A_{n;k}(L t_{0})
\left(\begin{array}{ll}
  {\alpha} & {\beta} \\
  {\gamma} & {\delta}
\end{array}\right)
\left(\begin{array}{l}
  {v(t)} \\
  {u(t)}
\end{array}\right) +
$$
$$
 + \xi_{L}(t) \sum_{\{j:k_{j}\neq 0\}} \int_{L t_{0}}^{L t} q_{j}(x)
\left(\begin{array}{ll}
  {v_{n;k-e_{j}}(x)} & {0} \\
  {u_{n;k-e_{j}}(x)} & {0}
\end{array}\right)
\left(\begin{array}{ll}
  {-u(x)} & {v(x)} \\
  {0} & {0}
\end{array}\right)
 d x
$$
$$
\left(\begin{array}{ll}
  {\alpha} & {\beta} \\
  {\gamma} & {\delta}
\end{array}\right)
\left(\begin{array}{l}
  {v(t)} \\
  {u(t)}
\end{array}\right) =
 \xi_{L}(t) A_{n;k}(L t_{0})
\left(\begin{array}{ll}
  {\alpha} & {\beta} \\
  {\gamma} & {\delta}
\end{array}\right)
\left(\begin{array}{l}
  {v(t)} \\
  {u(t)}
\end{array}\right) +
$$
$$
 + \xi_{L}(t) \sum_{\{j:k_{j}\neq 0\}} \int_{t_{0}}^{ t} q_{j}(s)
\frac{\xi_{L}(s)}{L'(s)} \left(\begin{array}{ll}
  {v_{n;k-e_{j}}(L s)} & {0} \\
  {u_{n;k-e_{j}}(L s)} & {0}
\end{array}\right)
\left(\begin{array}{ll}
  {-u(s)} & {v(s)} \\
  {0} & {0}
\end{array}\right)
 ds
\left(\begin{array}{l}
  {v(t)} \\
  {u(t)}
\end{array}\right)
$$
$$
 = \xi_{L}(t) A_{n;k}(L t_{0})
\left(\begin{array}{ll}
  {\alpha} & {\beta} \\
  {\gamma} & {\delta}
\end{array}\right)
\left(\begin{array}{l}
  {v(t)} \\
  {u(t)}
\end{array}\right) +
$$
$$
 + \xi_{L}(t) \sum_{\{j:k_{j}\neq 0\}}\int_{t_{0}}^{t} q_{j}(s)
[A_{n-1;k-e_{j}}(L t_{0}) \left(\begin{array}{ll}
  {\alpha} & {\beta} \\
  {\gamma} & {\delta}
\end{array}\right)
\left(\begin{array}{l}
  {v(s)} \\
  {u(s)}
\end{array}\right) +
$$
$$
 + \sum_{\{j_{1}:(k-e_{j})_{j_{1}}\neq 0\}}
A_{n-2;k-e_{j}-e_{j_{1}}}(L t_{0}) \left(\begin{array}{ll}
  {\alpha} & {\beta} \\
  {\gamma} & {\delta}
\end{array}\right)
\left(\begin{array}{l}
  {v_{1;e_{j_{1}}}(s)} \\
  {u_{1;e_{j_{1}}}(s)}
\end{array}\right) +
$$
$$
 + \sum_{\{j_{1}:(k-e_{j})_{j_{1}} \neq 0 \}}
\sum_{\{j_{2}:(k-e_{j}-e_{j_{1}})_{j_{2}} \neq 0 \}}
A_{n-3;k-e_{j}-e_{j_{1}}-e_{j_{2}}}(L t_{0})
\left(\begin{array}{ll}
  {\alpha} & {\beta} \\
  {\gamma} & {\delta}
\end{array}\right)
\left(\begin{array}{l}
  {v_{2;e_{j_{1}}+e_{j_{2}}}(s)} \\
  {u_{2;e_{j_{1}}+e_{j_{2}}}(s)}
\end{array}\right)
$$
$$
 + ... + \sum_{\{j_{1}:(k-e_{j})_{j_{1}}\neq 0\}}
\sum_{\{j_{2}:(k-e_{j}-e_{j_{1}})_{j_{2}} \neq 0 \}} ...
$$
$$
\sum_{\{j_{n-1}:(k-e_{j}-e_{j_{1}}-...-e_{j_{n-2}})_{j_{n-1}} \neq
0\}} A_{0;k-e_{j}-e_{j_{1}}-e_{j_{2}}-...-e_{j_{n-1}}}(L t_{0})
\left(\begin{array}{ll}
  {\alpha} & {\beta} \\
  {\gamma} & {\delta}
\end{array}\right)
$$
$$
\left(\begin{array}{l}
  {v_{n-1;e_{j_{1}}+e_{j_{2}}+...+e_{j_{n-1}}}(s)} \\
  {u_{n-1;e_{j_{1}}+e_{j_{2}}+...+e_{j_{n-1}}}(s)}
\end{array}\right) +
 \left(\begin{array}{ll}
  {\alpha} & {\beta} \\
  {\gamma} & {\delta}
\end{array}\right)
\left(\begin{array}{l}
  {v_{n;k-e_{j}}(s)} \\
  {u_{n;k-e_{j}}(s)}
\end{array}\right) ]
$$
$$
 \left(\begin{array}{ll}
  {-u(s)} & {v(s)} \\
  {0} & {0}
\end{array}\right)
d s \left(\begin{array}{l}
  {v(t)} \\
  {u(t)}
\end{array}\right) =
 \xi_{L}(t) A_{n;k}(L t_{0})
\left(\begin{array}{ll}
  {\alpha} & {\beta} \\
  {\gamma} & {\delta}
\end{array}\right)
\left(\begin{array}{l}
  {v(t)} \\
  {u(t)}
\end{array}\right) +
$$
$$
 + \xi_{L}(t) [\sum_{\{j:k_{j}\neq 0\}}
A_{n-1;k-e_{j}}(L t_{0}) \left(\begin{array}{ll}
  {\alpha} & {\beta} \\
  {\gamma} & {\delta}
\end{array}\right)
\left(\begin{array}{l}
  {v_{1;e_{j}}(t)} \\
  {u_{1;e_{j}}(t)}
\end{array}\right) +
$$
$$
 + \sum_{\{j:k_{j}\neq 0\}} \sum_{\{j_{1}:(k-e_{j})_{j_{1}}\neq 0\}}
A_{n-2;k-e_{j}-e_{j_{1}}}(L t_{0}) \left(\begin{array}{ll}
  {\alpha} & {\beta} \\
  {\gamma} & {\delta}
\end{array}\right)
\left(\begin{array}{l}
  {v_{2;e_{j}+e_{j_{1}}}(t)} \\
  {u_{2;e_{j}+e_{j_{1}}}(t)}
\end{array}\right) +
$$
$$
 + \sum_{\{j:k_{j}\neq 0\}} \sum_{\{j_{1}:(k-e_{j})_{j_{1}}\neq 0\}}
\sum_{\{j_{2}:(k-e_{j}-e_{j_{1}})_{j_{2}} \neq 0\}}
A_{n-3;k-e_{j}-e_{j_{1}}-e_{j_{2}}}(L t_{0})
\left(\begin{array}{ll}
  {\alpha} & {\beta} \\
  {\gamma} & {\delta}
\end{array}\right)
$$
$$
\left(\begin{array}{l}
  {v_{3;e_{j}+e_{j_{1}}+e_{j_{2}}}(t)} \\
  {u_{3;e_{j}+e_{j_{1}}+e_{j_{2}}}(t)}
\end{array}\right) + ... +
 \sum_{\{j:k_{j} \neq 0\}} \sum_{\{j_{1}:(k-e_{j})_{j_{1}} \neq 0\}}
\sum_{\{j_{2}:(k-e_{j}-e_{j_{1}})_{j_{2}} \neq 0\}} ...
$$
$$
\sum_{\{j_{n-1}:(k-e_{j}-e_{j_{1}}-...-e_{j_{n-2}})_{j_{n-1}} \neq
0\}} A_{0;k-e_{j}-e_{j_{1}}-e_{j_{2}}-...-e_{j_{n-1}}}(L t_{0})
\left(\begin{array}{ll}
  {\alpha} & {\beta} \\
  {\gamma} & {\delta}
\end{array}\right)
$$
$$
\left(\begin{array}{l}
  {v_{n;e_{j}+e_{j_{1}}+e_{j_{2}}+...+e_{j_{n-1}}}(t)} \\
  {u_{n;e_{j}+e_{j_{1}}+e_{j_{2}}+...+e_{j_{n-1}}}(t)}
\end{array}\right) +
 \left(\begin{array}{ll}
  {\alpha} & {\beta} \\
  {\gamma} & {\delta}
\end{array}\right)
\left(\begin{array}{l}
  {v_{n+1;k}(t)} \\
  {u_{n+1;k}(t)}
\end{array}\right) ].
$$
Taking into account the previous equality we  receive that the
following equality is valid
$$
\xi_{L}(t) \left(\begin{array}{ll}
  {\alpha_{L}(h)} & {\beta_{L}(h)} \\
  {\gamma_{L}(h)} & {\delta_{L}(h)}
\end{array}\right)
\left(\begin{array}{l}
  {v(t, h)} \\
  {u(t, h)}
\end{array}\right) =
\left(\begin{array}{l}
  {v(L t, h)} \\
  {u(L t, h)}
\end{array}\right) =
$$
$$
 = \left(\begin{array}{l}
  {v(L t)} \\
  {u(L t)}
\end{array}\right) +
\sum_{|k|=1} h^{k} \left(\begin{array}{l}
  {v_{1;k}(L t)} \\
  {u_{1;k}(L t)}
\end{array}\right) +
$$
$$
 + \sum_{|k|=2} h^{k}
\left(\begin{array}{l}
  {v_{2;k}(L t)} \\
  {u_{2;k}(L t)}
\end{array}\right) + ... +
 \sum_{|k|=n} h^{k}
\left(\begin{array}{l}
  {v_{n;k}(L t)} \\
  {u_{n;k}(L t)}
\end{array}\right) + ... =
$$
$$
 = \xi_{L}(t)[\left(\begin{array}{ll}
  {\alpha} & {\beta} \\
  {\gamma} & {\delta}
\end{array}\right)
\left(\begin{array}{l}
  {v(t)} \\
  {u(t)}
\end{array}\right) +
\sum_{|k|=1} h^{k} A_{0;k}(L t_{0}) \left(\begin{array}{ll}
  {\alpha} & {\beta} \\
  {\gamma} & {\delta}
\end{array}\right)
\left(\begin{array}{l}
  {v(t)} \\
  {u(t)}
\end{array}\right) +
$$
$$
 + \sum_{|k|=1} h^{k} \left(\begin{array}{ll}
  {\alpha} & {\beta} \\
  {\gamma} & {\delta}
\end{array}\right)
\left(\begin{array}{l}
  {v_{1;k}(t)} \\
  {u_{1;k}(t)}
\end{array}\right) +
 \sum_{|k|=2} h^{k} A_{1;k}(L t_{0})
\left(\begin{array}{ll}
  {\alpha} & {\beta} \\
  {\gamma} & {\delta}
\end{array}\right)
\left(\begin{array}{l}
  {v(t)} \\
  {u(t)}
\end{array}\right) +
$$
$$
 + \sum_{|k|=2} h^{k} \sum_{\{j:k_{j}\neq 0\}} A_{0;k-e_{j}}(L t_{0})
\left(\begin{array}{ll}
  {\alpha} & {\beta} \\
  {\gamma} & {\delta}
\end{array}\right)
\left(\begin{array}{l}
  {v_{1;e_{j}}(t)} \\
  {u_{1;e_{j}}(t)}
\end{array}\right) +
$$
$$
 + \sum_{|k|=2} h^{k} \left(\begin{array}{ll}
  {\alpha} & {\beta} \\
  {\gamma} & {\delta}
\end{array}\right)
\left(\begin{array}{l}
  {v_{2;k}(t)} \\
  {u_{2;k}(t)}
\end{array}\right) +
 \sum_{|k|=3} h^{k} A_{2;k}(L t_{0})
\left(\begin{array}{ll}
  {\alpha} & {\beta} \\
  {\gamma} & {\delta}
\end{array}\right)
\left(\begin{array}{l}
  {v(t)} \\
  {u(t)}
\end{array}\right) +
$$
$$
 + \sum_{|k|=3} h^{k} \sum_{\{j:k_{j}\neq 0\}} A_{1;k-e_{j}}(L t_{0})
\left(\begin{array}{ll}
  {\alpha} & {\beta} \\
  {\gamma} & {\delta}
\end{array}\right)
\left(\begin{array}{l}
  {v_{1;e_{j}}(t)} \\
  {u_{1;e_{j}}(t)}
\end{array}\right) +
$$
$$
 + \sum_{|k|=3} h^{k} \sum_{\{j:k_{j}\neq 0\}}
\sum_{\{(k-e_{j})_{j_{1}}\neq 0\}} A_{0;k-e_{j}-e_{j_{1}}}(L
t_{0}) \left(\begin{array}{ll}
  {\alpha} & {\beta} \\
  {\gamma} & {\delta}
\end{array}\right)
\left(\begin{array}{l}
  {v_{2;e_{j}+e_{j_{1}}}(t)} \\
  {u_{2;e_{j}+e_{j_{1}}}(t)}
\end{array}\right) +
$$
$$
 + \sum_{|k|=3} h^{k} \left(\begin{array}{ll}
  {\alpha} & {\beta} \\
  {\gamma} & {\delta}
\end{array}\right)
\left(\begin{array}{l}
  {v_{3;k}(t)} \\
  {u_{3;k}(t)}
\end{array}\right) + ... +
 \sum_{|k|=n} h^{k} A_{n-1;k}(L t_{0})
\left(\begin{array}{ll}
  {\alpha} & {\beta} \\
  {\gamma} & {\delta}
\end{array}\right)
\left(\begin{array}{l}
  {v(t)} \\
  {u(t)}
\end{array}\right) +
$$
$$
 + \sum_{|k|=n} h^{k} \sum_{\{j:k_{j}\neq 0\}} A_{n-2;k-e_{j}}(L t_{0})
\left(\begin{array}{ll}
  {\alpha} & {\beta} \\
  {\gamma} & {\delta}
\end{array}\right)
\left(\begin{array}{l}
  {v_{1;e_{j}}(t)} \\
  {u_{1;e_{j}}(t)}
\end{array}\right) +
$$
$$
 + \sum_{|k|=n} h^{k} \sum_{\{j:k_{j}\neq 0\}}
\sum_{\{j_{1}:(k-e_{j})_{j_{1}}\neq 0\}}
A_{n-3;k-e_{j}-e_{j_{1}}}(L t_{0}) \left(\begin{array}{ll}
  {\alpha} & {\beta} \\
  {\gamma} & {\delta}
\end{array}\right)
\left(\begin{array}{l}
  {v_{2;e_{j}+e_{j_{1}}}(t)} \\
  {u_{2;e_{j}+e_{j_{1}}}(t)}
\end{array}\right) + ... +
$$
$$
 \sum_{|k|=n} h^{k} \sum_{\{j:k_{j}\neq 0\}}
\sum_{\{j_{1}:(k-e_{j})_{j_{1}}\neq 0\}}
...\sum_{\{j_{n-2}:(k-e_{j}-e_{j_{1}}-...-e_{j_{n-3}})_{j_{n-2}}\neq
0\}} A_{0;k-e_{j}-e_{j_{1}}-...-e_{j_{n-2}}}(L t_{0})
$$
$$
\left(\begin{array}{ll}
  {\alpha} & {\beta} \\
  {\gamma} & {\delta}
\end{array}\right)
\left(\begin{array}{l}
  {v_{n-1;e_{j}+e_{j_{1}}+...+e_{j_{n-2}}}(t)} \\
  {u_{n-1;e_{j}+e_{j_{1}}+...+e_{j_{n-2}}}(t)}
\end{array}\right) +
 \sum_{|k|=n} h^{k}
\left(\begin{array}{ll}
  {\alpha} & {\beta} \\
  {\gamma} & {\delta}
\end{array}\right)
\left(\begin{array}{l}
  {v_{n;k}(t)} \\
  {u_{n;k}(t)}
\end{array}\right) + ... =
$$
$$
 = \xi_{L}(t)[
\left(\begin{array}{ll}
  {\alpha} & {\beta} \\
  {\gamma} & {\delta}
\end{array}\right)
\left(\begin{array}{l}
  {v(t, h)} \\
  {u(t, h)}
\end{array}\right) +
 \sum_{|k|=1} h^{k} A_{0;k}(L t_{0})
\left(\begin{array}{ll}
  {\alpha} & {\beta} \\
  {\gamma} & {\delta}
\end{array}\right)
\left(\begin{array}{l}
  {v(t, h)} \\
  {u(t, h)}
\end{array}\right) +
$$
$$
 \sum_{|k|=2} h^{k} A_{1;k}(L t_{0})
\left(\begin{array}{ll}
  {\alpha} & {\beta} \\
  {\gamma} & {\delta}
\end{array}\right)
\left(\begin{array}{l}
  {v(t, h)} \\
  {u(t, h)}
\end{array}\right) +
 \sum_{|k|=3} h^{k} A_{2;k}(L t_{0})
\left(\begin{array}{ll}
  {\alpha} & {\beta} \\
  {\gamma} & {\delta}
\end{array}\right)
\left(\begin{array}{l}
  {v(t, h)} \\
  {u(t, h)}
\end{array}\right)
$$
$$
 + ... +  \sum_{|k|=n} h^{k} A_{n-1;k}(L t_{0})
\left(\begin{array}{ll}
  {\alpha} & {\beta} \\
  {\gamma} & {\delta}
\end{array}\right)
\left(\begin{array}{l}
  {v(t, h)} \\
  {u(t, h)}
\end{array}\right) + ...].
$$

Thus, we received main result

{\bf Theorem 3.3.} {\sl{Let $\mathcal{M}[z(t, h)]$ be the
monodromy group for the solution  $z(t, h)$ of the Schwartz
equation (19) with normalization (5) on compact Riemann surface $F
= U/\Gamma.$ Then valid the exact variational formulas for $h
=(h_{1},..., h_{d}), |h| < \varepsilon,$
$$
\left(\begin{array}{ll}
  {\alpha_{L}(h)} & {\beta_{L}(h)} \\
  {\gamma_{L}(h)} & {\delta_{L}(h)}
\end{array}\right) =
[I + \sum _{|k|=1} A_{0;k}(Lt_{0})h^{k} +
 \sum _{|k|=2} A_{1;k}(Lt_{0}) h^{k} + ...
$$
$$
 + \sum _{|k|=n} A_{n-1;k}(Lt_{0}) h^{k} + ...]
\left(\begin{array}{ll}
  {\alpha_{L}(0)} & {\beta_{L}(0)} \\
  {\gamma_{L}(0)} & {\delta_{L}(0)}
\end{array}\right),
$$
where $I$ be unit matrix of order 2, $L \in \Gamma,$ $A_{0;k}(t) =
\int_{t_{0}}^{t}M_{j}(s) ds$ for $k = e_{j},$ $M_{j}(s) = q_{j}(s)
\left(\begin{array}{ll}
  {-u(s)v(s)} & {v(s)v(s)} \\
  {-u(s)u(s)} & {u(s)v(s)}
\end{array}\right),$
$j = 1,..., d;$
$$
A_{n;(k_{1},..., k_{d})}(t) = \int_{t_{0}}^{t} \sum_{\{j:k_{j}\neq
0\}}    A_{n-1;k-e_{j}}(s)M_{j}(s) ds,
$$
$n \geq 1, |k| \geq 2, d = 3g - 3, |h| < \varepsilon.$}}

Remark 3.2. The variational formulas show us how  the monodromy
group and the solution of the Schwartz equation depend of
accessory parameters $(h_{1},..., h_{d}).$ In particular, they
give the exact variational formula for generators any
quasifuchsian group and any Koebe group which uniformize a compact
Riemann surface of genus $g \geq 2.$

Remark 3.3. From a relation
$$
d A_{0}(L t) = \left(\begin{array}{ll}
  {\alpha_{L}(0)} & {\beta_{L}(0)} \\
  {\gamma_{L}(0)} & {\delta_{L}(0)}
\end{array}\right)
d A_{0}(t) \left(\begin{array}{ll}
  {\alpha_{L}(0)} & {\beta_{L}(0)} \\
  {\gamma_{L}(0)} & {\delta_{L}(0)}
\end{array}\right)^{-1},
L \in \Gamma, t \in U,
$$
follows that the matrix differential form $d A_{0}(t)$ is a Prym
differentials in sense of Gunning [25, с.160, 224], with respect
to characters, equal the monodromy homomorphism for Fuchsian group
$\Gamma$ with values in $PSL(2, {\bf C}).$ Thus, the variational
theory of monodromy groups for linearly polymorphic functions
(complex projective structures) on compact Riemann surface depends
on the periods such Prym differentials.

\centerline{{\bf References}}

[1] Hejhal D.A.  Monodromy groups and linearly polymorphic
function // Acta Math. 1975. V.135:1-2. P.1 - 55.

[2] Hejhal D.A. The variational theory of linearly polymorphic
functions // J. d'Analyse Math. 1976. V.30. P.215 - 264.

[3] Hejhal D.A. Kernel functions, Poincare series and LVA //
Contemporary Math. 2000. V.256. P.173 - 201.

[4] Earle C.J.  On variation of projective structures // Annals of
Math. Stud. N 97., Acad. Press., New-York. 1981. P.87 - 99.

[5] Kra I. Deformation of Fuchsian groups // Duke Math. J. 1969.
V.36. P.537 - 546.

[6] Kra I. Remarks on projective structures // Annals of Math.
Stud. N 97., Acad. Press., New-York. 1981. P.343 - 359.

[7] Maskit B. Uniformization of Riemann surface // Discontinuous
groups and Riemann surfaces, Ann. of  Math. Studies, Acad. Press.,
New-York. 1974. N 79. P.293 - 312.

[8] Gunning R.C.  Special coordinate coverings of Riemann surfaces
//  Math. Ann. 1967. V.170. P.67 - 86.

[9] Zograf P.G., Tachtajian L.A. Uniformization of Riemann
surfaces and the Weil-Petersson metric on Teichmueller and
Schottky spaces // Math. Sbornik. 1987. V.132, N 3. P.304 - 321.

[10] Venkov A.B. Examples effective solution of a Riemann-Hilbert
problem about restoration of the differential equations on to the
monodromy group within the framework of the theory automorphic
functions // Zapiski scientific seminars LOMI, S-Peterburg. 1987.
N 162. P.5 - 42

[11] Kapovich M.E. Hyperbolic manifolds and discrete groups.
Berlin: Progress in Mathematics, Birkhaeuser, 2001. V.183

[12] Appell P., Goursat E., Fatou P. Theorie des fonctions
algebriques. New-York: Chelsea Publish. Company, 1976.

[13] Ahlfors L.V., Bers L. Spaces of Riemann surfaces and
quasiconformal mapping. Moscow: IL, 1961.

[14] Maskit B. On the classification of Kleinian groups. II -
Signatures // Acta Math. 1977. V.138, N 1-2. P.17 - 42.

[15] Maskit B. Self-maps of Kleinian groups // Amer. J. Math.
1971. V. 93. P. 840 - 856.

[16] Forster O. Riemann surfaces. Moscow: Mir, 1980.

[17] Springer J. Introduction in the theory of Riemann surfaces.
Moscow: IL, 1960.

[18] Bers L. Holomorphic differentials as functions of moduli //
Bull. of the Amer. Math. Soc. 1961. V.67, N 2. P.206 - 210.

[19] Grauert H. Analytische Faserungen ueber
holomorphvollstaendigen Raeumen // Math. Ann. 1958. V.135. P.266 -
273.

[20] Krushkal S.L.  Quasiconformal mapping and Riemann surfaces
 Novosibirsk: Nauka, 1975.

[21] Maskit B. On the classification of Kleinian groups. I - Koebe
groups // Acta Math. 1975. V.135, N 3-4. P.249 - 271.

[22] Chueshev V.V. Spaces of compact Riemann surfaces and Koebe
groups // Sibirsk. Math. Jour. 1981 V.22, N 5. P.190 - 205.

[23] Golubev V.V. Lecture on the analytical theory of the
differential equations. Moscow - S.-Petersburg : GITTL, 1950.

[24] Shabat B.V. Introduction in the complex analysis, P.2.
Moscow: Nauka, 1985.

[25] Gunning R.C. Lectures on vector bundles over Riemann
surfaces. Princeton: Princeton Univ. Press., 1967.

\end{document}